\documentclass[11pt]{amsart}

\usepackage{etex}
\usepackage{amsmath,amssymb,caption,dbnsymb,enumerate,float,graphicx,kbordermatrix,multicol,multirow,setspace}
\usepackage[T1,T2A]{fontenc}
\usepackage[english]{babel}
\usepackage{pdfpages}
\usepackage{hyperref}\hypersetup{colorlinks,
	linkcolor={green!50!black},
	citecolor={green!50!black},
	urlcolor=blue}
\usepackage{epsfig}
\usepackage{mwe}
\usepackage{picins}
\usepackage{tikz}
\usepackage{tikz-cd}
\usetikzlibrary{arrows, matrix, positioning,calc}
\usepackage[all]{xy}

\newtheorem{theorem}{Theorem}[section]
\newtheorem{lemma}[theorem]{Lemma}
\newtheorem{proposition}[theorem]{Proposition}

\newtheorem{definition}[theorem]{Definition}
\newtheorem{example}[theorem]{Example}

\newtheorem{remark}[theorem]{Remark}

  \newcommand{\nc}{\newcommand}
\nc{\Xo}{X^{\text{out}}}
\nc{\Xin}{X^{\text{in}}}	
\nc{\al}{\alpha}
\nc{\be}{\beta}
\nc{\ga}{\gamma}
\nc{\de}{\delta}
\nc{\ep}{\epsilon}
\nc{\la}{\lambda}
\nc{\si}{\sigma}

\nc{\R}{\mathbb{R}}
\nc{\Z}{\mathbb{Z}}

\nc{\tMVA}{\text{tMVA}}
\nc{\rMVA}{\text{rMVA}}
\nc{\vMVA}{\text{vMVA}}
\nc{\nin}{n^{\text{in}}}
\nc{\no}{n^{\text{out}}}
\nc{\tr}{\text{tr}}
\nc{\xo}{x^{\text{out}}}
\nc{\xin}{x^{\text{in}}}

\newcommand {\DS} [1] {\displaystyle #1}

\begin{document}

\title{Alexander type invariants of tangles}

\author[Iva Halacheva]{Iva Halacheva}
\address{Department of Mathematics and Statistics \\
	Lancaster University \\ 
	Fylde College\\
	Lancaster, LA1 4YF, UK}
\email{i.halacheva@lancaster.ac.uk}
\urladdr{http://www.math.toronto.edu/ivahal}

\date{\today}

	\setlength{\parindent}{0pt}
	
	\begin{abstract} 
		We study generalizations of a classical link invariant -- the multivariable Alexander polynomial -- to tangles. The starting point is Archibald's tMVA invariant for virtual tangles which lives in the setting of circuit algebras, and whose target space has dimension that is exponential in the number of strands. Using the Hodge star map and restricting to tangles without closed components, we define a reduction of the tMVA to an invariant ``rMVA'' which is valued in matrices with Laurent polynomial entries, and so has a much more compact target space. We show the rMVA has the structure of a metamonoid morphism and is further equivalent to a tangle invariant defined by Bar-Natan. This invariant also reduces to the Gassner representation on braids and has a partially defined trace operation for closing open strands of a tangle.
	\end{abstract}
	
	\maketitle

	\tableofcontents

	\section{Introduction}
	
	A classical invariant of knots and links (i.e. multiple knotted copies of $S^1$) is the Alexander polynomial, originally introduced by James Alexander in 1928, \cite{Alex28}. A generalization of it to a multivariable polynomial invariant of links, with a different variable for each strand, was later presented by Torres \cite{Tor53} and is known as the multivariable Alexander polynomial, or MVA. Here, we will work with a version of it, vMVA (discussed more in Section \ref{sec:braids}), which is an invariant of oriented, regular, long, \textit{virtual} knots and links, i.e. those sitting in a thickened higher-genus surface rather than $\R^3$, first introduced by Kauffman \cite{K99}. A collection of objects which is more general than knots and links is that of {\sl tangles}, having knotted components which are either a smooth embedding of $S^1$ or of an interval, with its two endpoints fixed in two different planes. In her thesis \cite{Arch10}, Archibald defines an extension of the Alexander polynomial on links, called tMVA--an invariant of oriented, regular, virtual tangles from which can be recovered both the MVA and vMVA when closing the open strands of the tangle (compare also to \cite{B12}, \cite{DF16}, \cite{P10}). In fact, it satisfies the additional ``Overcrossing Commute'' relation pictured in Figure \ref{fig:R1OC}, which makes it an invariant of so-called welded tangles. This invariant is computed from the same type of Alexander matrix constructed for the MVA and provides a convenient setting for proving most of the local relations satisfied by the classical invariant it generalizes (see \cite{Arch10}). However, a drawback is that the dimension of its target space is exponential in the number of strands. More recently, Bar-Natan defines several versions of another, more compact, invariant for tangles, called $\beta$-calculus in \cite{BNS13} and \cite{DBN1}, and $\Gamma$-calculus in \cite{DBN2}. It originates in invariants of ribbon-knotted copies of $S^1$ and $S^2$ in $\mathbb{R}^4$ (see \cite{DBN3}), and associates to a pure tangle (i.e. one with no closed components), a (scalar, matrix) pair with Laurent polynomial entries. In this paper, we show that although the two invariants come from, on first sight, very different places, they carry essentially the same information for pure tangles. For the purpose, we show that the Archibald invariant is in fact determined by a small, much more compact part in the case of pure tangles. Namely, we start with Archibald's invariant in the setting of circuit algebras and use the Hodge star operator to transform it when restricted to pure tangles to a reduced, more compact version we will denote by rMVA. This reduced invariant fits in the algebraic structure of metamonoids, is readily computable, and also recovers the MVA. In addition, it is closer in shape to the invariant of Bar-Natan, and we produce in the context of metamonoids an isomorphism between their target spaces. The main result of this work is that, after the listed transformations, the two invariants are equivalent as metamonoid morphisms. The result provides a more efficient way of computing the tMVA, through the matrix-valued reduction rMVA. In addition, it provides a more direct proof that the Bar-Natan invariant is also a generalization of the Alexander polynomial. We further show these invariants restrict to the Gassner (and Burau) representations on braids, and have a partially defined trace operation allowing the closure of individual tangle components. The fact that they are readily computable by decomposing a tangle into its building blocks, and map into Laurent polynomials, also makes these invariants prime candidates for categorification.
	In Section \ref{sec:defres}, we discuss the preliminaries on tangles, circuit algebras, and the tMVA tangle invariant, as well as set the stage for reducing it. In Section \ref{sec:red} we define the reduced invariant rMVA and derive its values on positive and negative crossings, as well as the gluing, or multiplication, operations in the target space. In Section \ref{sec:meta}, we discuss the algebraic structure of metamonoids, show that the rMVA is a metamonoid morphism, recall the Bar-Natan invariant and show the two are equivalent as metamonoid morphisms. Recovering the original MVA invariant for links and the Gassner representation for braids from rMVA is described in Section \ref{sec:braids}. Some further results, including a partial trace operation and future directions, are discussed in Section \ref{sec:ext}. Some of the longer proofs are delegated to the end in Section \ref{sec:proof}.

	\subsection{Acknowledgements}
	This paper would not have been possible without the direction, support and encouragement of Dror Bar-Natan and I am extremely indebted to him. I am also grateful to the Knot at Lunch group at the University of Toronto, and in particular the multiple discussions with Ester Dalvit. This project was partially supported by an NSERC CGS-D scholarship.

	\section{Preliminaries}
\label{sec:defres}

\subsection{Circuit algebras and tangles}

The tangle invariant tMVA described by Archibald in \cite{Arch10} is defined in the context of circuit algebras, described in detail by Bar-Natan and Dancso \cite{BND2}, so we will start by discussing those. The operations in a circuit algebra are defined via circuit diagrams:

\begin{definition}[\cite{BND2}] An \textbf{oriented circuit diagram} (\textbf{OCD}) encodes the oriented pairing among a collection of ``input'' or ``internal'' and ``output'' or ``external'' points. 
	
	More precisely, let $\left[n_i\right]$ denote a set with $n_i$ elements, say $\{a_1,a_2,\hdots,a_{n_i}\}$, where $n_i \in \mathbb{Z}_{\geq 0}$. Given integers $k, \; l \in \mathbb{Z}_{\geq 0}$ and pairs $(n_i,m_i) \in \mathbb{Z}^2_{\geq 0}$, for $i=0, 1, \hdots, k$, an OCD prescribes a pairing among the elements of the set $\sqcup^k_{i=0}\left[n_i\right]\cup\left[m_i\right]$, and $l$ counts the number of closed oriented loops (we will elaborate on that later). The pairing must be oriented in the sense that a ``source'' point in $\sqcup^k_{i=0}\left[n_i\right]$ must be paired with a ``target'' point in $\sqcup^k_{i=0}\left[m_i\right]$ and conversely (so in particular if $\sum^k_{i=0}n_i \neq \sum^k_{i=0}m_i$ then no such pairing exists). Moreover, we think of $[n_0]\sqcup [m_0]$ as the ``output'' points and the rest, $\sqcup^k_{i=1}[n_i]\cup[m_i]$, as the ``input'' in the OCD. Formally, an OCD is an oriented compact $1$-manifold with boundary $(\left[n_i\right],\left[m_i\right]), i=0,1,\hdots,k$, up to homeomorphism.
	\vspace{-0.1cm}
	\parpic[r][b]{%
		\begin{minipage}{60mm}
			\begin{figure}[H]
				\begin{center}
					\def\svgwidth{4cm}
					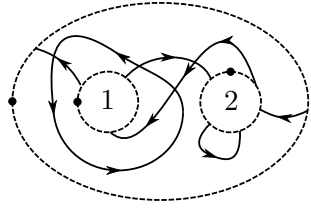
				\end{center}
				\vspace{0.2cm}
				\captionsetup{width=0.9\textwidth}
				\caption{An example of an oriented circuit diagram.}
				\label{fig:ocdeg}
			\end{figure}
		\end{minipage}
	}
	One way to represent such an OCD is via one ``external'' circle with marked points $(\left[n_0\right],\left[m_0\right])$ and $k$ ``internal'' circles, with marked points $(\left[n_i\right],\left[m_i\right]), i=1,\hdots,k$, together with $l$ oriented closed loops. In each such pair $(\left[n_i\right],\left[m_i\right])$, the first number indicates ``arrow tails'' and the second ``arrow heads''. We also number and place a dot on each circle from which we go counterclockwise when considering the marked points. Then, an OCD can be expressed as a pairing of the marked points by oriented arrows connecting them combinatorially (i.e. it is not important how they intersect) to the same circle or to another, which can only start at an arrow tail and end at an arrow head.  Informally, it is a pairing of the points $\sqcup^k_{i=0}\left[n_i\right]$ and $\sqcup^k_{i=0}\left[m_i\right]$, much like an electric circuit board with the ``internal'' circles being placeholders for chips. In Figure \ref{fig:ocdeg} can be seen an example of an OCD with $k=2,\; l=1,\; n_0=1,\; m_0=1,\; n_1=2,\; m_1=1,\; n_2=2,\; m_2=3$.
	
	Oriented OCDs can also be composed, as long as the marked points match. More precisely, we can compose an OCD $D$, whose parameters are given by $(k;l; \left(\left[n_i\right],\left[m_i\right]), i=0,1,\hdots,k\right)$, with OCDs $D_1,D_2,\hdots,D_k$ having parameters $(d_i;l_i; (\left[n_{ij}\right],\left[m_{ij}\right]), j=0,1,\hdots,d_i)$, $i=1,\hdots, k$, to obtain an OCD $D(D_1,D_2,\hdots,D_k)$ with parameters: \vspace{-0.2cm}$$\left(\sum^k_{i=1}d_i;\sum^k_{i=1}l_i+l+l';([n_0],[m_0]),(\left[n_{ij}\right],\left[m_{ij}\right]), i=1,\hdots, k, j=1,\hdots,d_i\right)$$\vspace{-0.2cm} 
	
	precisely when $n_i=m_{i0}, m_i=n_{i0} \; \forall \; i=1,\hdots,k$. Here $l' \in \mathbb{Z}_{\geq 0}$ denotes the number of new closed loops created from the composition. In Figure \ref{fig:ocds} we give an example of a composition of compatible OCDs.
	
	\begin{figure}[H]
		\begin{center}\def\svgwidth{12.3cm}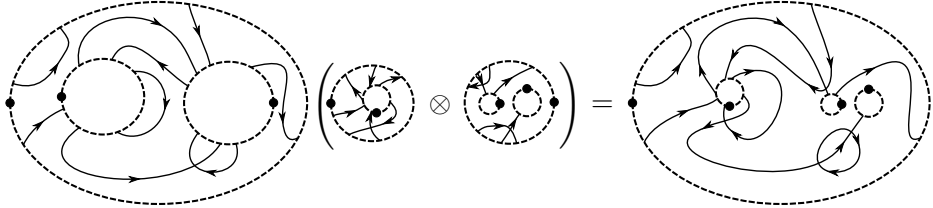\end{center}
		\caption{An example of a composition of oriented circuit diagrams.}
		\label{fig:ocds}
	\end{figure}
\end{definition}

The strands connect combinatorially to the same circle or to another, and the diagram represents an operation $V_{n_1,m_1} \otimes \hdots \otimes V_{n_k,m_k} \longrightarrow V_{n_0,m_0}$ (or $V_{n_1,m_1} \times \hdots \times V_{n_k,m_k} \longrightarrow V_{n_0,m_0}$ if we are dealing with sets), as we discuss below.

These combinatorial diagrams correspond to algebraic operations in the structure of an {\sl oriented circuit algebra}, which we define next.

\begin{figure}[H]
	\begin{center}\def\svgwidth{9cm}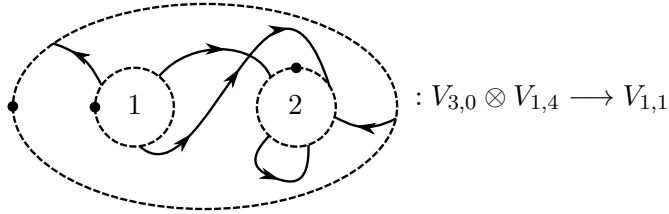\end{center}
	\caption{A circuit diagram representing an operation on a pair of objects.}
	\label{fig:ocd}
\end{figure}

\begin{definition}[\cite{BND2}]
	An \textbf{oriented circuit algebra} is an algebraic structure whose operations are indexed by oriented circuit diagrams. Namely, it is a collection $\mathcal{V}$ of objects and a collection $\mathcal{F}$ of operations, where:
	\begin{itemize}
		\item For each pair $(n,m)\in \mathbb{Z}^2_{\geq 0}$ there is a collection of objects $V_{n,m}$, which can have the structure of a set, vector space, module, etc.
		\item To every oriented circuit diagram $D$ (see Figure \ref{fig:ocd}), there is a corresponding operation, or morphism (of sets, vector spaces, modules, etc.) called $F_D$.
	\end{itemize}
	
\end{definition}

Analogously to planar algebras for usual tangles, introduced by Vaughan Jones \cite{Jon99}, circuit algebras provide a natural setting for discussing virtual tangles.

\begin{definition}[\cite{BND2}] Virtual tangle diagrams $v\mathcal{TD}$ have a presentation as the circuit algebra with two generators under strand concatenation: $\text{CA}\left\langle \overcrossing, \; \undercrossing \right\rangle$. Regular virtual or v-tangles $v\mathcal{T}$ are the circuit algebra quotient:
	$$\text{CA}\left\langle \overcrossing, \; \undercrossing \right\rangle/(\text{real Reidemeister moves})$$
	where the real Reidemeister moves for regular tangles are (with all possible orientations for the strands):
	\begin{figure}[H]\label{fig:reidemeister}	
		\begin{center}\def\svgwidth{10cm}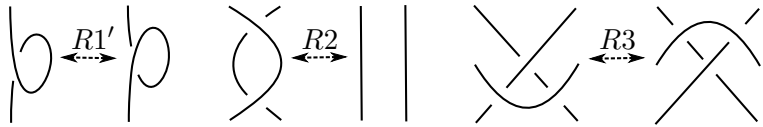\end{center}
		\caption{The real Reidemeister moves on diagrams for regular v-tangles.}
	\end{figure}
	In this oriented circuit algebra quotient:
	\begin{itemize}
		\item The collections of objects are indexed by $(n,m) \in \mathbb{Z}^2_{\geq 0}$ such that:
		\[V_{n,m}=\left\{
		\begin{array}{l}
		\emptyset, \text{ if } n \neq m \\
		\text{tangles with $n$ open strands and possibly} \\ 
		\text{some closed components}, \text{ if }  n=m 
		\end{array} 
		\right.\]
		\item The ``gluing'' operation within circuit diagrams is the concatenation of tangle strands.
		\begin{figure}[H]\label{fig:ocdtangles}	
			\begin{center}\def\svgwidth{12.5cm}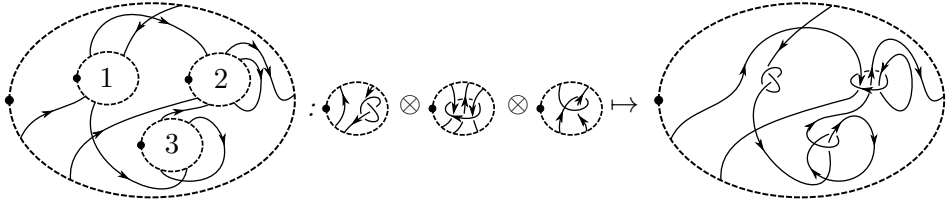\end{center}
			\caption{A circuit algebra operation on three tangles producing a composite tangle.}
		\end{figure}
	\end{itemize}
\end{definition}

The circuit algebra $v\mathcal{T}$ of virtual tangles is also graded by the skeleta (i.e. the underlying permutations) of the tangles.
\begin{remark}
	The generators of virtual tangle diagrams usually include an additional virtual crossing $\virtualcrossing$ and the virtual and mixed Reidemeister moves are imposed in order to represent virtual tangles (see Figure \ref{fig:vreid}). In the setting of circuit algebras, these are already built into the structure and come for free. We do not need to specifically include and require them.
	\begin{figure}[H]	
		\begin{center}\def\svgwidth{12.3cm}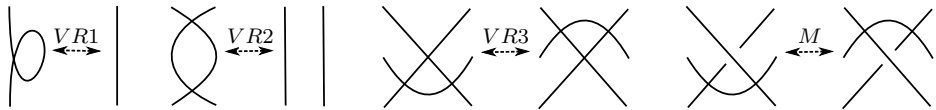\end{center}
		\caption{The virtual Reidemeister moves and the mixed move, which are part of the structure of circuit diagrams.}
		\label{fig:vreid}
	\end{figure}
	
	\parpic[r][b]{%
		\begin{minipage}{65mm}
			\begin{figure}[H]	
				\begin{center}\def\svgwidth{5cm}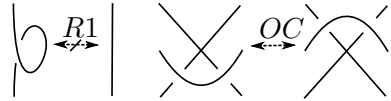\end{center}
				\captionsetup{width=0.9\textwidth}
				\caption{The usual R1 move does not hold for regular v-tangles, while the OC move gives w-tangles.}
				\label{fig:R1OC}
			\end{figure}
		\end{minipage}}
		Furthermore, note that the usual Reidemeister 1 relation does not hold in this setting. We will occasionally want to discuss the additional ``Overcrossings Commute'' (OC) relation (see Figure \ref{fig:R1OC}), and the quotient of virtual tangles by the OC relation is known as {\bf welded or w-tangles}, $w\mathcal{T}$.
	\end{remark}
	
	\subsection{The construction of tMVA}
	
	The target space for Archibald's tMVA invariant \cite{Arch10} is also a circuit algebra, $\mathcal{AHD}$, with the only nonempty sets of objects being $V_n:=V_{n,n}$, $n \in \mathbb{N}$, given by the {\bf Alexander half density} spaces:
	$$V_n=\text{AHD}([\nin],[\no]):=\Lambda^n([\no]) \otimes \Lambda^n([\nin] \cup [\no])$$
	where $[\nin]=\{a_1, \; a_2, \hdots, \; a_n\}$ and $[\no]=\{b_1, \; b_2, \hdots, b_n\}$ are $n$-element sets. Note that for $k \in \mathbb{Z}_{\geq 0}$ and $S$ a finite set, \vspace{-0.1cm}
	\begin{center}$\Lambda^k(S):=$ the $k$-th exterior power of the vector space with (formal) basis $S$.\footnote{Strictly speaking, for a tangle with $K$ strands in total (open and closed), we are taking the exterior power of the $\mathbb{R}[t_1,\hdots,t_K]$-module with basis $S$.}
	\end{center}
	
	The morphisms are given by ``gluing'' via interior multiplication. Namely, for a morphism $F_D: V_{n_1} \otimes \hdots \otimes V_{n_m} \longrightarrow V_{n_0}$ with gluing prescribed by a circuit diagram $D$, assign the same labels to the elements to be glued in $p_j \otimes q_j \in V_{n_j},\; j=1, \hdots, m$,\footnote{More formally, we also need to keep track of the closed components in the tangle.} denote that set $S=\bigcup^{m}_{j=1}[\nin_j] \cap \bigcup^{m}_{j=1}[\no_j]$.
	Then the gluing map is:
	$$i_S\left(\bigwedge^{m}_{j=1}p_j\right) \otimes i_S\left(\bigwedge^m_{j=1}q_j\right) \in \text{AHD}\left(\bigcup^{m}_{j=1}[\nin_j]-S, \bigcup^{m}_{j=1}[\no_j]-S\right)\vspace{0.1cm}$$
	where $i_S$ is interior multiplication with respect to $S$.
	
	The tMVA invariant then gives a circuit algebra morphism between v-tangles and Alexander half densities:
	$$\tMVA: (v\mathcal{T},\text{concatenation}) \longrightarrow (\mathcal{AHD}, \text{interior multiplication})$$
	
	To compute the tMVA on a regular, oriented v-tangle $T$, we first need the \textbf{Alexander matrix} $M(D_T)$ for a diagram $D_T$ of $T$. The matrix is indexed by the \textit{arcs} of the tangle, i.e. segments of the strands beginning and ending at an undercrossing or the outer circle. Let $\Xin$ be the set labelling the incoming arcs, and $\Xo$ the set of outgoing arcs of the tangle, $|\Xin|=|\Xo|=n$, and ``internal'' denote the remaining arc labels. Then the Alexander matrix is of the form:
	\begin{figure}[H]\vspace{-0.3cm}
		$$M(D_T)=\kbordermatrix{ &\text{internal} & & \Xo & & \Xin \cr
			&    & \vrule & & \vrule & \cr
			\text{internal} & & \vrule &  & \vrule & \cr	
			&  & \vrule & & \vrule & \cr \cline{2-6}
			&  & \vrule & & \vrule & \cr 
			\Xo &  & \vrule & & \vrule &} \hspace{1.5cm}$$
		\caption{The structure of the Alexander matrix $M(D_T)$ of a tangle diagram $D_T$.}
		\label{fig:strMT}
	\end{figure}

\begin{figure}[h]
		\begin{center}
			\begin{minipage}[c]{0.22\textwidth}
				\def\svgwidth{2.5cm}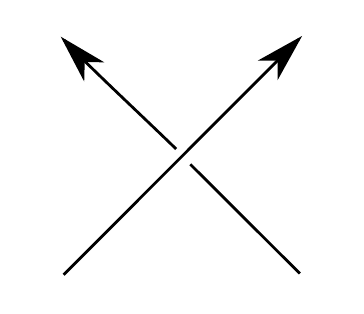
			\end{minipage}
			\begin{minipage}[c]{0.35\textwidth}
				$\mapsto \;
				\begin{array}{c|ccc}
				\phantom{s_0} & s_1 & s_2 & s_3 \\
				\hline
				s_3 & 1-t_2 & -1 & t_1
				\end{array}$
			\end{minipage}
			
			\vspace{0.5cm}
			
			\begin{minipage}[c]{0.22\textwidth}
				\def\svgwidth{2.5cm}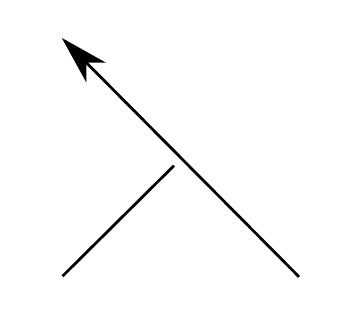
			\end{minipage}
			\begin{minipage}[c]{0.35\textwidth}
				$\mapsto \;
				\begin{array}{c|ccc}
				\phantom{s_0} & s_1 & s_2 & s_3 \\
				\hline
				s_3 & t_2-1 & -t_1 & 1
				\end{array}$
			\end{minipage}
		\end{center}
		
		\vspace{0.5cm}
		
		\begin{center}
			\begin{minipage}[c]{0.425\textwidth}
				\def\svgwidth{5cm}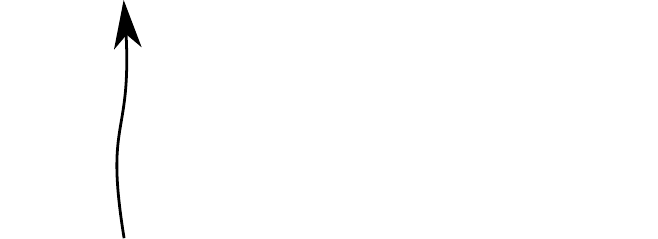
			\end{minipage}
			\begin{minipage}[c]{0.30\textwidth}
				$\mapsto \;
				\begin{array}{c|cc}
				\phantom{s_0} & a & b \\
				\hline
				a & 1 & -1
				\end{array}$
			\end{minipage}
		\end{center}
		\caption{The local rules for building the Alexander matrix $M(D_T)$ of a tangle diagram $D_T$.}
		\label{fig:rules}
	\end{figure}

	\begin{remark}
		To guarantee we have an equal number of distinct incoming and outgoing labels, we might also need to artificially break an arc $b$ into two differently labeled arcs $b$ and $a$. This results in an additional row and column of the matrix, both indexed by $a$, and does not affect the final expression for the invariant (see Figure \ref{fig:rules} and \cite{Arch10}, Lemma 3.9).
	\end{remark}
	
	In addition, we assign a variable to each strand of the tangle, so we have $\{t_i\}^{n+m}_{i=1}$ where $n$ is the number of open and $m$ is the number of closed components. The rows of the matrix are then obtained as below for each local piece of the tangle. If a column index does not appear, it is understood that the corresponding entry in the given row is zero.
	
	Here $s_i$ denotes the label of an arc, and $t_j$ is the variable corresponding to the $j^{\text{th}}$ strand of the tangle. From the Alexander matrix, Archibald defines the following tangle invariant:
	
	$$\boxed{\tMVA(T)= \prod^{n+m}_{s=1}{t^{-\frac{\mu(s)}{2}}_s} w \otimes \sum^n_{k=0} \hspace{-0.6cm}\sum_{\substack{\quad \quad \overline{i}=\{i_1 < \hdots < i_{n-k}\} \\
				\quad \quad \underline{j}=\{j_1 < \hdots < j_{k}\}\phantom{-k}}} \hspace{-0.6cm} \det{M(D_T)^{\overline{i};\underline{j}}} \hspace{0.2cm} b_{\overline{i}} \wedge a_{\underline{j}}\vspace{0.2cm}}$$ 
	
	The formula consists of the following ingredients:
	\begin{enumerate}[1)]
		\item $\mu(s)$ counts the number of times the $s^{\text{th}}$ strand is the overstrand when it passes through a crossing.
		\item $w \in \Lambda^n(\Xo)$ is determined by a choice of ordering of the elements of $\Xo$ as they appear in the rows (and columns) of the Alexander matrix.
		\item $M(D_T)^{\overline{i};\underline{j}}=M(D_T)^{i_1,\hdots,i_{n-k};j_1,\hdots,j_k}$ is the submatrix  of $M(D_T)$ with columns indexed by all the internal arcs, as well as $\{b_{i_1},\hdots,b_{i_{n-k}}\} \subset \Xo$ and $\{a_{j_1},\hdots,a_{j_k}\} \subset \Xin$.
		\item $b_{\overline{i}}=b_{i_1} \wedge \hdots \wedge b_{i_{n-k}}$ and $a_{\underline{j}}=a_{j_1} \wedge \hdots \wedge a_{j_k}$.
	\end{enumerate}
	
	The normalizing factor $\prod^{n+m}_{s=1}{t^{-\frac{\mu(s)}{2}}_s}$ belongs to $\mathbb{R}[t^{\pm\frac{1}{2}}_1,\hdots,t^{\pm\frac{1}{2}}_{n+m}]$, while the remaining expression lives in the Alexander half density space: \[\text{AHD}(\Xin,\Xo)=\Lambda^n(\Xo) \otimes \Lambda^{n}(\Xin \cup \Xo)\] 
	
	\begin{example} Below we compute the Alexander matrix for the given tangle diagram $D_T$ and the tMVA invariant of the corresponding tangle $T$. The expression might appear considerable but the computations are straightforward and can be delegated to a computer. 
		\vspace{0.2cm}
		\begin{figure}[H]
			\begin{center}
				\def\svgwidth{8cm}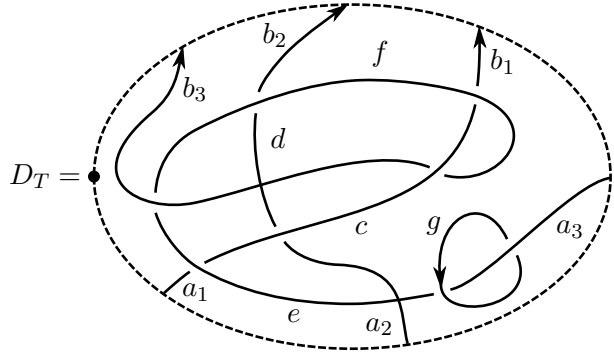
			\end{center}
			\caption{A tangle diagram $D_T$ for the tangle $T$.}
			\label{fig:tmvaeg1}
		\end{figure}
		\label{ex:tangle}
	\end{example}

	The corresponding Alexander matrix $M(D_T)$ for the tangle diagram $D_T$ is then given by:
				$$\begin{array}{c|ccccccccccc}
				&  c    &  d  &  e    &   f   &   g   & b_1 & b_2 &  b_3  &  a_1 &  a_2  & a_3 \\
				\hline
				c   & 1     &  0  & t_1-1 &   0   &   0   &  0  &  0  &   0   & -t_3 &  0 & 0 \\
				d   & 1-t_2 & t_1 &   0   &   0   &   0   &  0  &  0  &   0   &   0  & -1 & 0 \\
				e   & 0     & 0   &   1   &   0   & t_3-1 &  0  &  0  &   0   &   0  &  0 & -t_4 \\
				f   & 0     & 0   &  -t_3 &   1   &   0   &  0  &  0  & t_3-1 &   0  &  0 & 0 \\
				g   & 0     & 0   &   0   &   0   & 1-t_3 &  0  &  0  &   0   &   0  &  0 & t_4-1 \\
				b_1 & -1    & 0   &   0   & 1-t_1 &   0   & t_3 &  0  &   0   &   0  &  0 &  0  \\
				b_2 & 0     & -1  &   0   & 1-t_2 &   0   &  0  & t_3 &   0   &   0  &  0 &  0  \\
				b_3 & 1-t_3 & 0   &   0   &  -1   &   0   &  0  &  0  &  t_1  &   0  &  0 &  0 
				\end{array}$$
				
	From the Alexander matrix, together with the normalizing term for the tangle diagram $D_T$, we get the value tMVA(T) of the tMVA invariant:	
	
								$$\begin{array}{rl}
								&t^{-1}_1t^{-\frac{5}{2}}_3t^{-\frac{1}{2}}_4 \; b_1 \wedge b_2 \wedge b_3 \otimes (t_3-1)t_3\left[\right.(t_1-1)^2 (t_3-1) \; b_3\wedge a_2\wedge a_3+ \\
								&+(t_3 t_1-t_1-2 t_3+1)\; b_1\wedge a_2\wedge a_3-t_3 (t_3 t_1-t_1-t_3) \; a_1\wedge a_2\wedge a_3+\\
								&+(t_2-1)(t_3 t_1-t_1-1) t_3^2 \; b_1\wedge a_1\wedge a_3	+t_1 (t_3-1) t_3^2 \; b_1\wedge b_2\wedge a_1-\\
								&t_1 (t_3 t_1-t_1-t_3) t_3^2 \; b_2\wedge a_1\wedge a_3-t_1 (t_3 t_1-2 t_1-t_3+1) t_3^2 \; b_2\wedge b_3\wedge a_1+\\
								&+(t_3-1) t_3 \; b_1\wedge a_1\wedge a_2-t_1 (t_3 t_1-t_1-2 t_3+1) t_3 \; b_1\wedge b_2\wedge a_3+\\
								&+(t_1-1)^2 t_1 (t_3-1) t_3 \; b_2\wedge b_3\wedge a_3-t_1 (t_1+t_3-1) t_3 \; b_1\wedge b_2\wedge b_3+\\
								&+(t_3 t_1-2 t_1-t_3+1) t_3 \; b_3\wedge a_1\wedge a_2-t_1 (t_2-1) t_3 \; b_3\wedge a_1\wedge a_3-\\
								&-(t_1+t_3-1) \; b_1\wedge b_3\wedge a_2+(t_2-1)(t_1 t_3^2-2 t_1 t_3-t_3+1) t_3 \; b_1\wedge b_3\wedge a_1 \\
								&-(t_2-1)(t_3^2 t_1^2-t_3 t_1^2-t_3^2 t_1+t_3 t_1+t_1+t_3-1) \; b_1\wedge b_3\wedge a_3 \left.\right]
								\end{array}$$

	\subsection{The Hodge star operator}
	
	We can apply the Hodge star operator $\ast_w$ to the second tensor component of $\text{AHD}(\Xo,\Xin)$, for a given choice of $w \in  \Lambda^n(\Xo)$ (coming from a choice of ordering of the elements of $\Xo$), as follows: 
	
	
	\begin{equation}
	\begin{tikzcd}[ampersand replacement=\&]
	\DS{\bigoplus^n_{k=0}}  \Lambda^k(\Xin) \otimes \Lambda^{n-k}(\Xo)\arrow{rr}{\ast_w} \& \&
	\DS{\bigoplus^n_{k=0}} \Lambda^k(\Xin) \otimes \Lambda^k(\Xo) \arrow{d}{\cong} \\
	\Lambda^n(\Xin \cup \Xo)  \arrow[swap]{u}{\cong} \& \&
	\DS{\bigoplus^n_{k=0}} \text{Hom}(\Lambda^k(\Xin), \Lambda^k(\Xo)) \\
	\end{tikzcd}
	\label{eqn:Hodge}
	\end{equation}
	
	It acts as the identity on $\Lambda^k(\Xin)$, on the basis $\{b_{i_1}\wedge \hdots \wedge b_{i_{n-k}}\}_{1 \leq i_1<\hdots<i_{n-k} \leq n}$  of $\Lambda^{n-k}(\Xo)$ (where $\Xo=\{b_1,\hdots, b_n\}$) as described below, and we extend linearly ($i_{n-k+1} < \hdots < i_n$):
		\[\begin{array}{c}
		b_{i_1} \wedge \hdots \wedge b_{i_{n-k}} \mapsto (-1)^{\ast_w} b_{i_{n-k+1}} \wedge \hdots \wedge b_{i_n} \\
		\\
		 \Leftrightarrow \; (-1)^{\ast_w}b_{i_1} \wedge \hdots \wedge b_{i_{n-k}} \wedge b_{i_{n-k+1}} \wedge \hdots \wedge b_{i_n} = w
		\end{array}\]
		
	Let $p_2: \Lambda^n(\Xo) \otimes \Lambda^{n}(\Xin \cup \Xo) \longrightarrow \Lambda^{n}(\Xin \cup \Xo)$ denote the projection onto the second tensor factor in $\text{AHD}(\Xin,\Xo)$.
	
	If we consider, for a tangle $T$ and a given $w \in \Lambda^n(\Xo)$, the image of $p_2(\tMVA(T))$ in the space $\bigoplus^n_{k=0} \text{Hom}(\Lambda^k(\Xin), \Lambda^k(\Xo))$ under $\ast_w$, it turns out that this image is determined by the degree $0$ and $1$ components.
	
	\begin{theorem}
		For a tangle $T$, let $\lambda \in 
		\text{Hom}(\Lambda^0(\Xin), \Lambda^0(\Xo))$ denote the degree $0$ component, and $\phi\in 
		\text{Hom}(\Lambda^1(\Xin), \Lambda^1(\Xo))$ the degree $1$ component of the image of $p_2(\tMVA(T))$ in $\bigoplus^n_{k=0} \text{Hom}(\Lambda^k(\Xin), \Lambda^k(\Xo))$ under the Hodge star operator $\ast_w$. If $\lambda \neq 0$, then the whole image is determined by $\la$ and $\phi$ -- it is precisely $\lambda \cdot \Lambda(\phi/\lambda)$.
		\label{thm:deg01maps}
	\end{theorem}
	
	The proof will be postponed until Section \ref{sec:proof}.
	
	Let us now consider the degree $0$ and $1$ components of the image of $T$, $p_2(\tMVA(T))$, in the space $\bigoplus^n_{k=0} \Lambda^k(\Xin) \otimes \Lambda^k(\Xo)$ using $\ast_{w}$.

	The degree $0$ component is $\la = \det{M(D_T)^{1,\hdots,n;\emptyset}}$, where (recalling earlier notation) $M(D_T)^{1,\hdots,n;\emptyset}$ is the submatrix of $M(D_T)$ with columns indexed by all internal arcs, as well as all $n$ elements of $\Xo$ but none from $\Xin$.
	
	The degree $1$ elements can be expressed as the entries of an $n \times n$ matrix $\mathcal{A}$, where $\mathcal{A}_{i,j}$ is the determinant of $M(D_T)^{1,\hdots,n;\emptyset}$ with column $i$ in $\Xo$ replaced by column $j$ from $\Xin$:
	\[\mathcal{A}_{i,j}=(-1)^{n-i}\det{M(D_T)^{1,\hdots,\hat{i},\hdots,n;j}}\]
	This $n\times n$ matrix, with the additional degree $0$ element and the normalizing factor $\prod^n_{s=1}{t^{-\frac{\mu(s)}{2}}_s}$ in the definition of the tMVA invariant carries the information of $\tMVA(T)$.
	
	\begin{theorem}\label{thm:deg01}
		Let $T$ be a tangle with associated pair $(\lambda,\mathcal{A})$ as above. When $\lambda \neq 0$, we can recover from $(\lambda,\mathcal{A})$ the coefficients of the tangle invariant $\tMVA(T)$ via the formula:
		\[
		\det{M(D_T)^{\{1,\hdots,n\} \setminus \{i_1,\hdots,i_k\};j_1,\hdots,j_k}}=(-1)^{nk-\frac{(k-1)k}{2}-\sum^k_{p=1}{i_p}}\frac{\det{\mathcal{A}^{j_1,\hdots,j_k}_{i_1,\hdots,i_k}}}{\lambda^{k-1}}
		\]
		where $\mathcal{A}^{j_1,\hdots,j_k}_{i_1,\hdots,i_k}$ denotes the submatrix of $\mathcal{A}$ with rows $i_1,\hdots,i_k$ and columns $j_1,\hdots,j_k$.
	\end{theorem}

	This proof will also be postponed until Section \ref{sec:proof}. In the next section, we will explore this version of the tMVA invariant further.

	\section{The Reduction of tMVA to rMVA}
\label{sec:red}
Since the new pair $(\lambda,\mathcal{A})$ together with the normalizing factor preserves the information of tMVA, as expected it is also a tangle invariant that we'll denote rMVA (``r'' for reduced). Its computation can be simplified similarly to that of tMVA by breaking up a tangle into the generating pieces, i.e. positive and negative crossings, computing the invariant on each piece, and then gluing the pieces back together. For this purpose, we need to know the values of rMVA on positive and negative crossings, and its behaviour under gluing strands and taking the disjoint union of tangles. 

To define rMVA, we restrict to \textbf{pure} regular v-tangles, i.e. ones without closed components. Then in particular, as we'll see in Proposition \ref{prop:nzero}, $\la \neq 0$ and so Theorem \ref{thm:deg01} lets us recover tMVA from rMVA. Given such an $n$-tangle $T$ with incoming and outgoing labels for the strands being the sets $\Xin$ and $\Xo$ respectively, we can identify $\Xin \cong \Xo =X$ from the underlying permutation of the tangle, using the strand variables $\{t_i\}^n_{i=1}$. The set $X$ then also provides a labelling for the set of strands of the tangle. We will denote the set of pure regular v-tangles labelled using $X$ by $pv\mathcal{T}_X$.

If $(\la,\mathcal{A})=(\la_T, \mathcal{A}_T) \in \mathbb{R}(t_i) \times M_{X \times X}(\mathbb{R}(t_i))$ are the pair of degree $0$ and $1$ components for the image of a pure, regular v-tangle $T$ as defined above for a fixed $w \in \Lambda^n(\Xo)$, then Theorem \ref{thm:deg01} reduces tMVA to:
\begin{equation}\boxed{\rMVA(T):=\prod^n_{k=1}{t^{-\frac{\mu(k)}{2}}_k}(\la,\mathcal{A}) \in R_X:=\mathbb{R}(\sqrt{t_i}) \times M_{X \times X}(\mathbb{R}(\sqrt{t_i}))}\label{eqn:dirdef}\end{equation}
\subsection{A tangle invariant}

Since the tMVA is a virtual tangle invariant and rMVA is a function of it, it also has that property. We will nevertheless check for illustration purposes that it satisfies the real Reidemeister 2 and 3 moves, as well as the ``Overcrossings Commute'' relation, which makes it a welded tangle invariant. We include here the last verification and delegate the remaining ones to Section \ref{sec:proof}. We need to verify that the values of rMVA on each side of the move agree. Note that as the virtual crossings do not contribute to the Alexander matrix, rMVA is automatically invariant under the virtual and mixed Reidemeister moves.
\vspace{0.5cm}

\underline{\sl Overcrossings Commute moves.}
\vspace{-0.2cm}
\begin{figure}[H]
	\begin{center}	\hspace{-1.25cm}
		\begin{minipage}[c]{0.30\textwidth}
			\def\svgwidth{5cm}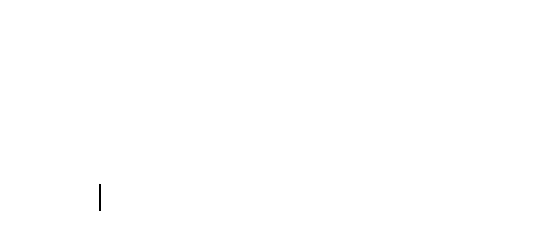
		\end{minipage}
		\hspace{2.5cm}
		\begin{minipage}[c]{0.30\textwidth}
			\def\svgwidth{5cm}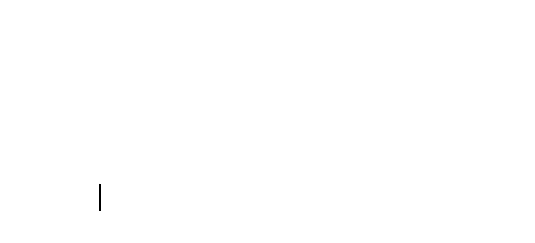
		\end{minipage}
	\end{center}
	\caption{The braid-like and cyclic Overcrossings Commute moves.}
	\label{fig:ocbc}
\end{figure}
\vspace{-0.5cm}
\begin{remark}
	We distinguish two types of Overcrossings Commute or OC moves in Figure \ref{fig:ocbc}, depending on whether the region bounded by the three strands can be oriented consistently from the strand orientations or not. If not, the move is called ``braid-like'', otherwise it is ``cyclic''.
\end{remark}
Starting with the braid-like move, the diagram on the left side $D_{LB}$ has the following Alexander matrix and corresponding value of the rMVA:
\vspace{-0.1cm}
\begin{align*}
	M_{LB} &=	\begin{array}{c|cccccc}
		& b_1 & b_2 & b_3 & a_1 & a_2 & a_3 \\
		\hline
		b_1 & 1 & 0 & t_1-1 & -t_3 & 0 & 0 \\
		b_2 & 0 & 1 & 0 & 0 & -t_3 & t_2-1 \\
		b_3 & 0 & 0 & 1 & 0 & 0 & -1
	\end{array} \\
	\Rightarrow
	\rMVA(D_{LB})&=\prod^3_{s=1}{t^{-\frac{\mu(s)}{2}}_s}\cdot
	\begin{array}{c|c}
		\la_{LB} & \Xin \\
		\hline
		\Xo & \mathcal{A}_{LB}
	\end{array}=
	t^{-1}_3 \cdot
	\begin{array}{c|ccc}
		1 & a_1 & a_2 & a_3\\
		\hline
		b_1 & -t_3 & 0 & t_1-1 \\
		b_2 & 0 & -t_3 & t_2-1 \\
		b_3 & 0 & 0 & -1
	\end{array}
\end{align*}
\vspace{-0.1cm}
The diagram on the right side of the braid-like Overcrossings Commute move, $D_{RB}$, produces:
\vspace{-0.1cm}
\begin{align*}
	M_{RB} &=	\begin{array}{c|cccccc}
		& b_1 & b_2 & b_3 & a_1 & a_2 & a_3 \\
		\hline
		b_1 & 1 & 0 & 0 & -t_3 & 0 & t_1-1 \\
		b_2 & 0 & 1 & t_2-1 & 0 & -t_3 & 0 \\
		b_3 & 0 & 0 & 1 & 0 & 0 & -1
	\end{array} \\
	\Rightarrow
	\rMVA(D_{RB})&=\prod^3_{s=1}{t^{-\frac{\mu(s)}{2}}_s}\cdot
	\begin{array}{c|c}
		\la_{RB} & \Xin \\
		\hline
		\Xo & \mathcal{A}_{RB}
	\end{array}=
	t^{-1}_3 \cdot
	\begin{array}{c|ccc}
		1 & a_1 & a_2 & a_3\\
		\hline
		b_1 & -t_3 & 0 & t_1-1 \\
		b_2 & 0 & -t_3 & t_2-1 \\
		b_3 & 0 & 0 & -1
	\end{array}
\end{align*}

Analogously, for the left side of the cyclic Overcrossings Commute move, we have:
\begin{align*}
	M_{LC} &=	\begin{array}{c|cccccc}
		& b_1 & b_2 & b_3 & a_1 & a_2 & a_3 \\
		\hline
		b_1 & 1 & 0 & t_1-1 & -t_3 & 0 & 0 \\
		b_2 & 0 & t_3 & 0 & 0 & -1 & 1-t_2 \\
		b_3 & 0 & 0 & 1 & 0 & 0 & -1
	\end{array} \quad 	\Rightarrow \\
	\rMVA(D_{LC})&=\prod^3_{s=1}{t^{-\frac{\mu(s)}{2}}_s}\cdot
	\begin{array}{c|c}
		\la_{LC} & \Xin \\
		\hline
		\Xo & \mathcal{A}_{LC}
	\end{array}=
	t^{-1}_3 \cdot
	\begin{array}{c|ccc}
		t_3 & a_1 & a_2 & a_3\\
		\hline
		b_1 & -t^2_3 & 0 & t_3(1-t_1) \\
		b_2 & 0 & -1 & 1-t_2 \\
		b_3 & 0 & 0 & -t_3
	\end{array}
\end{align*}

The diagram on the right side of the cyclic move has the same rMVA value:
\begin{align*}
	M_{RC} &=	\begin{array}{c|cccccc}
		& b_1 & b_2 & b_3 & a_1 & a_2 & a_3 \\
		\hline
		b_1 & 1 & 0 & 0 & -t_3 & 0 & t_1-1 \\
		b_2 & 0 & t_3 & 1-t_2 & 0 & -1 & 0 \\
		b_3 & 0 & 0 & 1 & 0 & 0 & -1
	\end{array} \quad \Rightarrow \\
	\rMVA(D_{RC})&=\prod^3_{s=1}{t^{-\frac{\mu(s)}{2}}_s}\cdot
	\begin{array}{c|c}
		\la_{RC} & \Xin \\
		\hline
		\Xo & \mathcal{A}_{RC}
	\end{array}=
	t^{-1}_3 \cdot
	\begin{array}{c|ccc}
		t_3 & a_1 & a_2 & a_3\\
		\hline
		b_1 & -t^2_3 & 0 & t_3(1-t_1) \\
		b_2 & 0 & -1 & 1-t_2 \\
		b_3 & 0 & 0 & -t_3
	\end{array}
\end{align*}

\begin{remark}
	Note that rMVA is a regular v-tangle invariant like tMVA, and does not satisfy the real Reidemeister 1 move in Figure \ref{fig:R1OC}. Indeed, comparing the values of rMVA on the two diagrams, we see that they differ:
	\begin{figure}[h]
		\begin{center}
			\begin{minipage}[c]{0.12\textwidth}
				\def\svgwidth{1.5cm}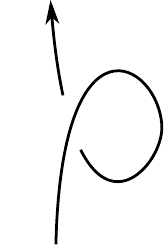
			\end{minipage}
			\begin{minipage}[c]{0.3\textwidth}
				\begin{align*}&\xrightarrow{M} 
					\begin{array}{c|cc}
						\phantom{a} & b_1 & a_1 \\
						\hline
						b_1 & t_1 & -t_1
					\end{array} \\
					&\xrightarrow{\rMVA}
					t^{-1/2}_1 \cdot \begin{array}{c|cc}
						t_1 & a_1 \\
						\hline
						b_1 & -t_1
					\end{array}\end{align*}
				\end{minipage}
				\hspace{1cm}
				\begin{minipage}[c]{0.03\textwidth}
					\def\svgwidth{0.55cm}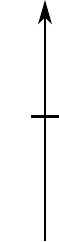
				\end{minipage}
				\begin{minipage}[c]{0.3\textwidth}
					\begin{align*}&\xrightarrow{M} 
						\begin{array}{c|cc}
							\phantom{a} & b_1 & a_1 \\
							\hline
							b_1 & 1 & -1
						\end{array} \\
						&\xrightarrow{\rMVA}
						1 \cdot \begin{array}{c|cc}
							1 & a_1 \\
							\hline
							b_1 & -1
						\end{array}\end{align*}
					\end{minipage}
				\end{center}
			\end{figure}
		\end{remark}
		
		After verifying that rMVA is indeed a tangle invariant, we simplify its description by recovering its values on positive and negative crossings. We also describe the strand gluing and disjoint union operations in the target space, induced from those for the tMVA invariant in the Alexander half density spaces.
		
		\subsection{Positive and negative crossings}
		
		The Alexander matrices for the positive and negative crossings are obtained after splitting the overcrossing arc into two. From them, we can find the value of rMVA:
		\begin{figure}[H]
			\begin{center}
				\begin{minipage}[c]{0.15\textwidth}
					\def\svgwidth{2cm}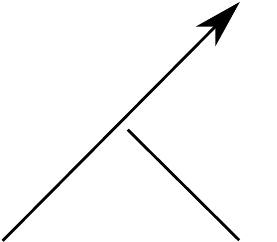
				\end{minipage}
				\begin{minipage}[c]{0.8\textwidth}
					$\xrightarrow{M} 
					\begin{array}{c|cccc}
					\phantom{a} & b_1 & b_2 & a_1 & a_2 \\
					\hline
					b_1 & 1 & 0 & -1 & 0 \\
					b_2 &  1-t_2 & t_1 & 0 & -1
					\end{array} \xrightarrow{\rMVA}
					t^{-1/2}_1 \cdot \begin{array}{c|cc}
					t_1 & a_1 & a_2 \\
					\hline
					b_1 & -t_1 & 0 \\
					b_2 & 1-t_2  & -1
					\end{array}
					$
				\end{minipage}
			\end{center}
			
			\vspace{0.5cm}
			
			\begin{center}
				\begin{minipage}[c]{0.15\textwidth}
					\def\svgwidth{2cm}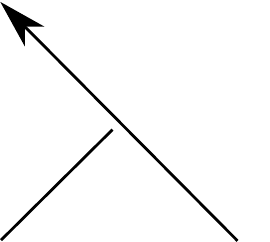
				\end{minipage}
				\begin{minipage}[c]{0.82\textwidth}
					$\xrightarrow{M}
					\begin{array}{c|cccc}
					\phantom{s_0} & b_1 & b_2 & a_1 & a_2 \\
					\hline
					b_1 & 1 & 0 & -1 & 0 \\
					b_2 & t_2-1 & 1 & 0 & -t_1
					\end{array} \xrightarrow{\rMVA}
					t^{-1/2}_1 \cdot\begin{array}{c|cc}
					1 & a_1 & a_2 \\
					\hline
					b_1 & -1 & 0 \\
					b_2 & t_2-1 & -t_1 
					\end{array}$
				\end{minipage}
			\end{center}
			\caption{The Alexander matrices and resulting degree $0$ and $1$ pairs for the positive and negative crossing tangles.}
			\label{fig:pnmat}
		\end{figure}

		\subsection{Gluing and disjoint union}
		\label{sub:glue}
		Let $T$ denote a pure tangle with labels $X\; (=\Xo\cong \Xin)$. The operations of gluing and disjoint union on $\tMVA(T)$ in the target space $\text{AHD}(X)$ of tMVA correspond to gluing two strands of the tangle $T$ or taking the disjoint union of two tangles, respectively. They can also be studied on the level of the Alexander matrix. We will do so to reduce them to analogous operations on the pairs below, in the target space $R_X=\mathbb{R}(\sqrt{t_i}) \times M_{X \times X}(\mathbb{R}(\sqrt{t}_i))$ of rMVA:
		\[\left\{\rMVA(T)=\prod^{|X|}_{s=1}{t^{-\frac{\mu(s)}{2}}_s}(\lambda_T,\mathcal{A}_T)\right\}_{T\; \in \; pv\mathcal{T}_X}\]
		\vspace{-0.1cm}
		For this purpose, we consider the map $M(D_T) \mapsto (\lambda_T,\mathcal{A}_T)$ taking the Alexander matrix to such a pair.
		\begin{figure}[H]
			\begin{center}
				\def\svgwidth{10cm}
				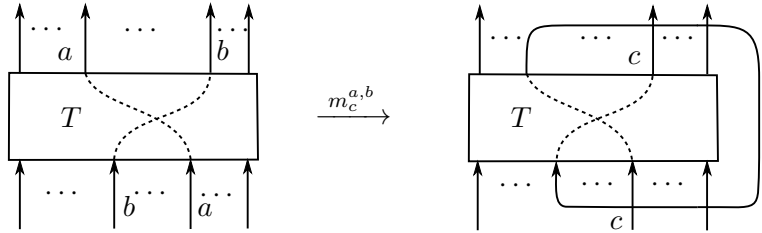
				\caption{Gluing the strands labelled ``$a$'' and ``$b$'' of a tangle $T$, and calling the resulting strand ``$c$''.}
				\label{fig:glue}
			\end{center}
		\end{figure}
		\vspace{-0.4cm}
		\begin{lemma}\textbf{(Gluing)}
			The result of gluing the outgoing strand labelled ``$a$'' $\in X$ in a tangle $T$ with ends labelled by $X$, to the incoming strand labelled ``$b$'' $\in X$, and calling the resulting strand ``$c$'', as illustrated in Figure \ref{fig:glue}, corresponds to the following map on the level of the target space $R_X$ of rMVA. (For the moment we restrict to gluing different strands so we do not get closed components.)
			\begin{figure}[H]
				\begin{equation}\arraycolsep=4pt
				\begin{array}{c|ccc}
				\lambda & a & b & \Xin \\
				\hline
				a & \alpha & \beta & \theta \\
				b & \gamma & \delta & \epsilon \\
				\Xo & \phi & \psi & \Xi 
				\end{array}
				\xrightarrow{\; m^{a,b}_c \;}
				\left(\arraycolsep=8pt\def\arraystretch{1.5}
				\begin{array}{c|ccc}
				\lambda + \beta & c & \Xin \\
				\hline
				c & \gamma + \frac{\beta\gamma - \alpha\delta}{\lambda} & \epsilon + \frac{\beta\epsilon - \delta\theta}{\lambda} \\
				\Xo & \phi + \frac{\beta\phi-\alpha\psi}{\lambda} & \Xi + \frac{\beta\Xi - \psi\theta}{\lambda} 
				\end{array}\right)_{t_a,t_b \rightarrow t_c}
				\label{eqn:glue}
				\end{equation}
			\end{figure}
			\label{lem:glue}
		\end{lemma}
		\vspace{-0.4cm}
		\begin{proof}
			Note that for the operation $m^{a,b}_c$ on $R_X$ in Equation \ref{eqn:glue}, any pair $(\la,\mathcal{A}) \in R_X$ as above, and any scalar $\mu$, we have $m^{a,b}_c(\mu\cdot(\la,\mathcal{A}))=(\mu)_{t_a,t_b \rightarrow t_c}\cdot m^{a,b}_c(\la,\mathcal{A})$, which will be of use when dealing with the normalizing factor. The gluing operation $m^{a,b}_c$ on a tangle $T$, pictured in Figure \ref{fig:glue}, has the following effect on the Alexander matrix, mapping $M(D_T)$ to $N(D_T)=m^{a,b}_c(M(D_T))$, where also $t_a,t_b \rightarrow t_c$:
			\vspace{-0.1cm}
			\[\begin{array}{c}
			\kbordermatrix{ &\overset{\text{\tiny internal}}{\phantom{a}} & & \overset{\normalsize \Xo}{a^{\text{out}} \; b^{\text{out}}} & & \overset{\normalsize \Xin}{a^{\text{in}} \; b^{\text{in}}} \cr
				\rotatebox{90}{\text{\tiny internal}}  & & \vrule &  & \vrule & \cr	\cline{2-6}		
				a^{\text{out}} &  & \vrule & & \vrule & \cr
				b^{\text{out}} &  & \vrule & & \vrule & \cr
				\Xo	&  & \vrule & & \vrule & }
			\xrightarrow{m^{a,b}_c} \\
			\xrightarrow{m^{a,b}_c}
			\kbordermatrix{ & \overset{\text{\tiny internal \quad }}{  l=a^{\text{out}}+b^{\text{in}}} & & \overset{\normalsize \Xo}{\widehat{a^{\text{out}}} \; c^{\text{out}}(=b^{\text{out}})} & & \overset{\normalsize \Xin}{c^{\text{in}}(=a^{\text{in}}) \; \widehat{b^{\text{in}}}} \cr
				\rotatebox{90}{\text{\tiny internal}} & & \vrule &  & \vrule & \cr	
				l=a^{\text{out}} &  & \vrule &  & \vrule &  \cr \cline{2-6} 
				\widehat{a^{\text{out}}} &  & \vrule & & \vrule &  \cr
				c^{\text{out}}(=b^{\text{out}}) &  & \vrule & & \vrule &  \cr
				\Xo &  & \vrule & & \vrule &}\end{array}\vspace{0.5cm}\]
			
			We denote by $M(D_T)$ and $\rMVA(T)=\prod^{|X|}_{s=1}{t^{-\frac{\mu(s)}{2}}_s}(\lambda,\mathcal{A})$ the original Alexander matrix and rMVA invariant values, and denote the result after gluing by $N(D_T):=m^{a,b}_c(M(D_T))$ for the new Alexander matrix, and the invariant: 
			\begin{align*}
			\rMVA(m^{a,b}_c(T))&=m^{a,b}_c\left(\prod^{|X|}_{s=1}{t^{-\frac{\mu(s)}{2}}_s}\cdot(\lambda,\mathcal{A})\right) \\ &=\left(\prod^{|X|}_{s=1}{t^{-\frac{\mu(s)}{2}}_s}\right)_{t_a,t_b\rightarrow t_c}\hspace{-0.9cm}m^{a,b}_c(\lambda,\mathcal{A})=:\left(\prod^{|X|}_{s=1}{t^{-\frac{\mu(s)}{2}}_s}\right)_{t_a,t_b\rightarrow t_c}\hspace{-1cm}(\omega,\mathcal{B})\end{align*}
			So, it suffices to verify Formula \ref{eqn:glue} for $(\omega, \mathcal{B})$. To obtain $N(D_T)$, we delete column $a^{\text{out}}$ in $\Xo$ and column $b^{\text{in}}$ in $\Xin$ and include a new internal arc with column labelled $l$, which is the sum of the original columns $a^{\text{out}}$ and $b^{\text{in}}$. Furthermore, we need to remove row $a^{\text{out}}$ in $\Xo$ and include it, now labeled $l$, with the rows labeled by the internal arcs. We denote the overall new strand obtained after gluing $a$ and $b$ by $c$. So, in the columns of the resulting matrix we have $c^{\text{out}}$ equal the old $b^{\text{out}}$ and $c^{\text{in}}$ being the old $a^{\text{in}}$. In the rows of the new matrix, $c^{\text{out}}$ is equal to the old $b^{\text{out}}$. Using the matrix $N(D_T)$, we can compute directly the resulting pair $(\omega, \mathcal{B})$ from the definition:
			\vspace{-0.6cm}
			\begin{figure}[H]
				\begin{align*}
					\omega &= \det{N(D_T)^{1,\hdots,n;\emptyset}} =\det{\kbordermatrix{ & \overset{\text{\tiny internal \quad }}{\phantom{a} l=a^{\text{out}}+b^{\text{in}}} & & \overset{\normalsize \Xo}{\widehat{a^{\text{out}}} \; c^{\text{out}}(=b^{\text{out}})} &  \cr
							\rotatebox{90}{\text{\tiny internal}} & & \vrule &  &  \cr	
							l=a^{\text{out}} &  & \vrule &  &   \cr \cline{2-4} 
							\widehat{a^{\text{out}}} &  & \vrule & &   \cr
							c^{\text{out}}(=b^{\text{out}}) &  & \vrule & &  \cr
							\Xo &  & \vrule & & }} \\
					& \\
					&=\det{\kbordermatrix{ & \overset{\text{\tiny internal \quad }}{\phantom{a^{\text{out}}} \quad \widehat{l}} & & \overset{\normalsize \Xo}{a^{\text{out}} \; b^{\text{out}}}  \cr
							\rotatebox{90}{\text{\tiny internal}} & & \vrule &    \cr	
							\widehat{l} &  & \vrule &     \cr \cline{2-4} 
							a^{\text{out}} &  & \vrule &    \cr
							b^{\text{out}} &  & \vrule &  \cr
							\Xo &  & \vrule &  }}  +
					\det{\kbordermatrix{ & \overset{\text{\tiny internal \quad }}{\phantom{a^{\text{out}}} \quad \widehat{l}} & & \overset{\normalsize \Xo}{b^{\text{in}} \; b^{\text{out}}}   \cr
							\rotatebox{90}{\text{\tiny internal}} & & \vrule &  \cr	
							\widehat{l} &  & \vrule &     \cr \cline{2-4} 
							a^{\text{out}} &  & \vrule &    \cr
							b^{\text{out}} &  & \vrule &  \cr
							\Xo &  & \vrule &  }}  \\
					& \\
					&=\la + \mathcal{A}_{a,b}=\la + \beta
				\end{align*}
			\end{figure}
			\vspace{-0.2cm}
			Similarly, we can obtain the entries of $\mathcal{B}$. By definition, we have that for $d \in (\Xo \setminus \{a\})_{b \rightarrow c}$ and $e \in (\Xin \setminus \{b\})_{a \rightarrow c}$, $\mathcal{B}_{d,e}$ is the determinant of the matrix $N(D_T)^{1,\hdots,n;\emptyset}$ with column $d$ in $\Xo$ replaced with column $e$ from $\Xin$:
			\vspace{-0.6cm}
			\begin{figure}[H]
				\begin{align*}
					&\mathcal{B}_{d,e} =\det{\kbordermatrix{ & \overset{\text{\tiny internal \quad }}{\phantom{a^{\text{out}}} l=a^{\text{out}}+b^{\text{in}}} & & \overset{\normalsize \Xo \phantom{c^{\text{out}}(=b^{\text{out}})} \quad \large{\widehat{d^{\text{out}}}}}{\widehat{a^{\text{out}}} \; c^{\text{out}}(=b^{\text{out}}) \; e^{\text{in}}} &  \cr
							\rotatebox{90}{\text{\tiny internal}} & & \vrule &  &  \cr	
							l=a^{\text{out}} &  & \vrule &  &   \cr \cline{2-5} 
							\widehat{a^{\text{out}}} &  & \vrule & &   \cr
							c^{\text{out}}(=b^{\text{out}}) &  & \vrule & &  \cr
							\Xo &  & \vrule & & }} \\
					\end{align*}
				\end{figure}

				\vspace{5cm}
				\begin{figure}[H]
				\begin{align*}
					&=\det{\hspace{-0.2cm}\kbordermatrix{ & \overset{\text{\tiny internal }}{\widehat{l}\hspace{-1cm}} &  & \overset{\normalsize \Xo \phantom{c} \quad \large{\widehat{d^{\text{out}}}}}{a^{\text{out}} \; b^{\text{out}} \; e^{\text{in}}}  \cr
							\rotatebox{90}{\text{\tiny internal}} & & \vrule &   \cr	
							\widehat{l} &  & \vrule &    \cr \cline{2-4} 
							a^{\text{out}} &  & \vrule &   \cr
							b^{\text{out}} &  & \vrule &  \cr
							\Xo &  & \vrule & }} +
					\det{\hspace{-0.2cm}\kbordermatrix{ & \overset{\text{\tiny internal }}{\widehat{l}\hspace{-1cm}} &  & \overset{\normalsize \Xo \phantom{c} \quad \large{\widehat{d}}}{b^{\text{in}} \; b^{\text{out}} \; e^{\text{in}}}  \cr
							\rotatebox{90}{\text{\tiny internal}} & & \vrule &  \cr	
							\widehat{l} &  & \vrule &   \cr \cline{2-4} 
							a^{\text{out}} &  & \vrule &   \cr
							b^{\text{out}} &  & \vrule &  \cr
							\Xo &  & \vrule & }}\\
					&= \mathcal{A}_{d',e'} + \frac{\det{\mathcal{A}^{b,e'}_{a,d'}}}{\la}
				\end{align*}
			\end{figure}
			
			where the last equality follows from the definition of the degree $1$ matrix $\mathcal{A}$ and Theorem \ref{thm:deg01}:
			\[
			d'=\left\{
			\begin{array}{ll}
			d, \text{ if } d \in \Xo \setminus \{a,b\} \\ 
			b, \text{ if } d=c
			\end{array}
			\right.
			\quad
			e'=\left\{
			\begin{array}{ll}
			e, \text{ if } e \in \Xin \setminus \{a,b\} \\ 
			a, \text{ if } e=c
			\end{array}
			\right.
			\]
			This agrees with the formula for the gluing map $m^{a,b}_c$ on $R_X$ given in Equation \ref{eqn:glue}.
		\end{proof}
		\vspace{0.4cm}
		\begin{lemma}\textbf{(Disjoint Union)}
			When taking the disjoint union of two tangles $T_1$ and $T_2$, i.e. putting them side by side to create a new tangle, the resulting map on the target space of rMVA is:
			\begin{equation}\arraycolsep=1pt\def\arraystretch{1.5}
			\begin{array}{c|c}
			\lambda_1 & \; \Xin_1 \; \\
			\hline
			\; \Xo_1 \; & \; A_1 \; \\
			\end{array}
			\cup
			\begin{array}{c|c}
			\lambda_2 & \; \Xin_2 \; \\
			\hline
			\; \Xo_2 \; & \; A_2 \; \\
			\end{array}
			=
			\begin{array}{c|cc}
			\lambda_1 \cdot\lambda_2 & \; \Xin_1 \; & \; \Xin_2 \; \\
			\hline
			\Xo_1 & \lambda_2 A_1 & 0 \\
			\Xo_2 & 0 &  \lambda_1 A_2 \\
			\end{array}\label{eqn:union}\end{equation}
			\label{lem:union}
		\end{lemma}
	
		\begin{proof}
			The Alexander matrix for a diagram of the disjoint union $T_1 \cup T_2$ of the two tangles will be of the form:
			
			\[\arraycolsep=5pt\def\arraystretch{1.5}
			M(D_{T_1 \cup T_2}) = \kbordermatrix{ & \text{int}_1 & & \text{int}_2 & & \Xo_1 & & \Xo_2 & & \Xin_1 & & \Xin_2 \cr
				\text{int}_1 &  & \vrule & O & \vrule & & \vrule & O & \vrule & & \vrule & O &  \cr	\cline{2-13}
				\text{int}_2 & O & \vrule &  & \vrule & O & \vrule & & \vrule & O & \vrule & &  \cr	\cline{2-13}
				\Xo_1 &  & \vrule & O & \vrule & & \vrule & O & \vrule & & \vrule & O &  \cr	\cline{2-13}
				\Xo_2 & O  & \vrule & & \vrule & O & \vrule &  & \vrule & O & \vrule & &}
			\]
			
			\vspace{0.4cm}
			
			Like the gluing operation, the disjoint union operation on $R_X \times R_X$ respects scalar multiplication: for any two pairs $(\la_1,\mathcal{A}_1), \; (\la_2,\mathcal{A}_2)  \in R_X$ and any scalars $\mu_1,\; \mu_2$, we have:
			\[\left[\mu_1(\la_1,\mathcal{A}_1)\right] \sqcup \left[\mu_2(\la_2,\mathcal{A}_2)\right]=\mu_1\mu_2\left[(\la_1,\mathcal{A}_1)\sqcup(\la_2,\mathcal{A}_2)\right]\]
			Suppose $|\Xo_1|=|\Xin_1|=n_1,\; |\Xo_2|=|\Xin_2|=n_2$, and $\rMVA(T_1)=\eta_1(\lambda_1,\mathcal{A}_1), \; \rMVA(T_2)=\eta_2(\lambda_2,\mathcal{A}_2)$, where $\eta_i$ are the normalizing coefficients and $\la_i, \; \mathcal{A}_i$ are the degree $0$ and $1$ components  respectively, under the Hodge star map. Then for the resulting invariant $\rMVA(T_1\cup T_2)=\rMVA(T_1)\cup\rMVA(T_2)=\eta_1\eta_2(\omega, \mathcal{B})$, we have:
			\[\omega = \arraycolsep=2pt\def\arraystretch{1.5}\det{\kbordermatrix{ & \text{int}_1 & & \text{int}_2 & & \Xo_1 & & \Xo_2  \cr
					\text{int}_1 &  & \vrule & O & \vrule & & \vrule & O &   \cr	\cline{2-9}
					\text{int}_2 & O & \vrule &  & \vrule & O & \vrule & &   \cr	\cline{2-9}
					\Xo_1 &  & \vrule & O & \vrule & & \vrule & O &  \cr	\cline{2-9}
					\Xo_2 & O  & \vrule & & \vrule & O & \vrule &  & }  	} 
			= \det{\kbordermatrix{ & \text{int}_1 & & \Xo_1 & & \text{int}_2 & & \Xo_2  \cr
					\text{int}_1 &  & \vrule &  & \vrule & O & \vrule & O &   \cr	\cline{2-9}
					\Xo_1 &  & \vrule &  & \vrule & O & \vrule & O &  \cr	\cline{2-9}
					\text{int}_2 & O & \vrule & O & \vrule & & \vrule & &   \cr	\cline{2-9}
					\Xo_2 & O  & \vrule & O & \vrule &  & \vrule &  & }}
			\]
			
			\vspace{0.5cm}
			
			Namely, $\omega=\lambda_1\lambda_2$. Furthermore, the $(i,j)$-th entry of the $(n_1+n_2) \times (n_1+n_2)$ matrix $\mathcal{B}$ is obtained by taking the determinant of $M(D_{T_1 \cup T_2})^{1,\hdots,{n_1+n_2};\emptyset}$ with the $i$-th column among the labels $\Xo_1 \cup \Xo_2$ replaced with the $j$-th column among the labels $\Xin_1 \cup \Xin_2$. If we replace the $i$-th column from $\Xo_2$ with the $j$-th column $a^1_j$ from $\Xin_1$, we get a zero determinant, i.e. $\mathcal{B}_{n_1+i,j}=0, \; \forall \; 1 \leq i \leq n_2,\; 1 \leq j \leq n_1$:
			
			\begin{align*}
				\mathcal{B}_{n_1+i,j}=&\det{\arraycolsep=3pt\def\arraycolstretch{1.3}\kbordermatrix{ & \text{int}_1 & & \text{int}_2 & & \Xo_1 &  & \Xo_2 & a^1_j & \cr
						\text{int}_1 &  & \vrule & O & \vrule & & \vrule & 0 & \ast & 0 \cr	\cline{2-10}
						\text{int}_2 & O & \vrule &  & \vrule & O & \vrule & \ast & 0 & \ast  \cr	\cline{2-10}
						\Xo_1 &  & \vrule & O & \vrule & & \vrule & 0 & \ast & 0 \cr	\cline{2-10}
						\Xo_2 & O  & \vrule & & \vrule & O & \vrule & \ast & 0 & \ast } } \\
				= &\det{\kbordermatrix{ & \text{int}_1 & & \Xo_1 & & \text{int}_2 & & \Xo_2 & a^1_j & \cr
						\text{int}_1 &  & \vrule &  & \vrule & O & \vrule & 0 & \ast & 0  \cr	\cline{2-10}
						\Xo_1 &  & \vrule &  & \vrule & O & \vrule & 0 & \ast & 0  \cr	\cline{2-10}
						\text{int}_2 & O & \vrule & O & \vrule & & \vrule & \ast & 0 & \ast   \cr	\cline{2-10}
						\Xo_2 & O  & \vrule & O & \vrule &  & \vrule &  \ast & 0 & \ast }} \\
				=&\det{\kbordermatrix{ & \text{int}_1 & & \Xo_1 \cr
						\text{int}_1 &  & \vrule &  \cr	\cline{2-4}
						\Xo_1 &  & \vrule &  \cr }}
				\det{\kbordermatrix{ & \text{int}_2 & & \Xo_2 & a^1_j & \cr
						\text{int}_2 & & \vrule & \ast & 0 & \ast   \cr	\cline{2-6}
						\Xo_2 &  & \vrule &  \ast & 0 & \ast }}=0
			\end{align*}
			
			\vspace{0.4cm}
			
		Analogously for a choice of columns from $\Xo_1$ and $\Xin_2$, we have that $\mathcal{B}_{i,n_1+j}=0, \; \forall \; 1 \leq i \leq n_1, \; 1 \leq j \leq n_2$. 
			
			In the more interesting case, if we replace the $i$-th column from $\Xo_1$ with the $j$-th column $a^1_j$ from $\Xin_1$, we get:
			
			\begin{align*}
				\mathcal{B}_{i,j}&=\det{\arraycolsep=3pt\def\arraycolstretch{1.3}\kbordermatrix{ & \text{int}_1 & & \text{int}_2 & & \Xo_1 & a^1_j &  & \Xo_2  \cr
						\text{int}_1 &  & \vrule & O & \vrule & &  &   \vrule & O \cr	\cline{2-9}
						\text{int}_2 & O & \vrule &  & \vrule & O & 0 &   \vrule &   \cr	\cline{2-9}
						\Xo_1 &  & \vrule & O & \vrule & &  &  \vrule & O \cr	\cline{2-9}
						\Xo_2 & O  & \vrule & & \vrule & O &  0 &  \vrule &  } } \\
				&= \det{\kbordermatrix{ & \text{int}_1 & & \Xo_1 & a^1_j & & \text{int}_2 & & \Xo_2  \cr
						\text{int}_1 &  & \vrule & & & \vrule & O &   \vrule & O \cr	\cline{2-9}
						\Xo_1 &  & \vrule & & & \vrule & O &  \vrule & O \cr	\cline{2-9}
						\text{int}_2 & O & \vrule & O & 0 & \vrule &  &   \vrule &   \cr	\cline{2-9}
						\Xo_2 & O & \vrule & O & 0 & \vrule &  &  \vrule &  } } \\
				&= \det{\kbordermatrix{ & \text{int}_1 & & \Xo_1 & a^1_j  \cr
						\text{int}_1 &  & \vrule & &  \cr	\cline{2-5}
						\Xo_1 &  & \vrule & & \cr} }
				\det{\kbordermatrix{& \text{int}_2 & & \Xo_2  \cr
						\text{int}_2 & & \vrule &    \cr	\cline{2-4}
						\Xo_2 & & \vrule &  } }	\\
				&= (\mathcal{A}_1)_{i,j}\cdot\lambda_2
			\end{align*}
			
			\vspace{0.2cm}
			
			Similarly, replacing the $i$-th column from $\Xo_2$ with the $j$-th column from $\Xin_2$, we get: \[\mathcal{B}_{n_1+i,n_1+j}=(\mathcal{A}_2)_{i,j}\cdot \lambda_1, \; \forall \; 1 \leq i,j \leq n_2\]
			\vspace{0.1cm}
			This agrees with the formula for disjoint union in Equation \ref{eqn:union}.
		\end{proof}
		
		\vspace{0.2cm}
		
		\begin{remark}
			For the gluing operation of rMVA to be well-defined and for this invariant to determine tMVA, we need to have $\la_T \neq 0$ when computing the normalized pair $(\la_T,\mathcal{A}_T)$ for a tangle $T$. Note that in the positive and negative crossings, if we set all the variables $t_i=1$, we get the pair $(1,-I_{2\times 2})$. Furthermore, the operations of gluing (as long as we are gluing different strands and not creating closed components) and disjoint union preserve the property $\la \neq 0$, as we show in Proposition \ref{prop:nzero}. So, for pure tangles it is always true that $\la \neq 0$.
		\end{remark}
		
		\begin{proposition}
			\label{prop:nzero}
			For any pure regular v-tangle $T$ with corresponding normalized degree $0$ and $1$ pair $(\la_T, \mathcal{A}_T)$, we always have $\la_T \neq 0$.
		\end{proposition}
		\begin{proof}
			It suffices to check that this property is true for the positive and negative crossing generators, and that it is preserved by the disjoint union and gluing maps. We have:
			\begin{figure}[H]
				\begin{center}
					\begin{minipage}{0.12\textwidth}
						\def\svgwidth{1.5cm}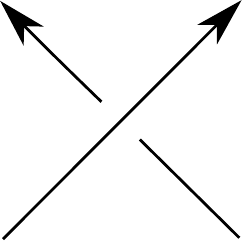
					\end{minipage}
					\begin{minipage}{0.7\textwidth}
						$\xrightarrow{\rMVA}
						t^{-1/2}_a \cdot \begin{array}{c|cc}
						t_a & a & b \\
						\hline
						a & -t_a & 0 \\
						b & 1-t_b  & -1
						\end{array}
						\xrightarrow{t_i=1}
						\begin{array}{c|cc}
						1 & a & b \\
						\hline
						a & -1 & 0 \\
						b & 0  & -1
						\end{array}
						$
					\end{minipage}
				\end{center}
			\end{figure}
			\vspace{-0.5cm}
			\begin{figure}[H]
				\begin{center}
					\begin{minipage}{0.12\textwidth}
						\def\svgwidth{1.5cm}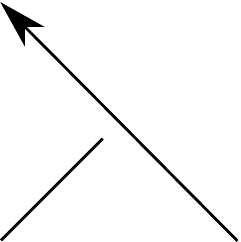
					\end{minipage}
					\begin{minipage}{0.70\textwidth}
						$\xrightarrow{\rMVA}
						t^{-1/2}_a \cdot\begin{array}{c|cc}
						1 & a & b \\
						\hline
						a & -1 & 0 \\
						b & t_b-1 & -t_a 
						\end{array}
						\xrightarrow{t_i=1}
						\begin{array}{c|cc}
						1 & a & b \\
						\hline
						a & -1 & 0 \\
						b & 0  & -1
						\end{array}
						$
					\end{minipage}
				\end{center}
			\end{figure}	
			In particular, since the $\la_T$ values for the positive and negative crossings satisfy $\la_{\overcrossing}(1)=\la_{\undercrossing}(1)=1$, they are not identically zero. We will show that for any $n$-component pure v-tangle $T$, $(\la_T,\mathcal{A}_T)_{t_i=1}=(1,-I_{n \times n})$. This will mean in particular that $\la_T$ is not identically zero. It suffices to check that this property is preserved under gluing and disjoint union.
			\begin{figure}[H]
				\[\arraycolsep=5pt
				\begin{array}{c|ccc}
				\lambda & a & b & X \\
				\hline
				a & \alpha & \beta & \theta \\
				b & \gamma & \delta & \epsilon \\
				X & \phi & \psi & \Xi 
				\end{array}
				\xrightarrow{\; m^{a,b}_c \;}
				\left(\arraycolsep=10pt\def\arraystretch{1.5}
				\begin{array}{c|ccc}
				\lambda + \beta & c & X \\
				\hline
				c & \gamma + \frac{\beta\gamma - \alpha\delta}{\lambda} & \epsilon + \frac{\beta\epsilon - \delta\theta}{\lambda} \\
				X & \phi + \frac{\beta\phi-\alpha\psi}{\lambda} & \Xi + \frac{\beta\Xi - \psi\theta}{\lambda} 
				\end{array}\right)_{t_a,t_b \rightarrow t_c}
				\]
			\end{figure}
			Assuming that in the original matrix $\la(1)=1, \; \alpha(1)=\delta(1)=-1, \; \Xi(1)=-I_{(n-2) \times (n-2)}$, $\beta(1)=\gamma(1)=0, \; \theta(1)=\epsilon(1)=(\phi(1))^{\text{tr}}=(\psi(1))^{\text{tr}}=(0,\hdots,0)$, the element and matrix pair resulting from the gluing operation $m^{a,b}_c$ evaluated at $t_i=1$, is of the form $(1,-I_{(n-1) \times (n-1)})$:
			
			\begin{align*}
				&(\la+\beta)(1) =1+0=1, \quad \left(\gamma+\frac{\beta\gamma - \alpha\delta}{\la}\right)(1)= 0 + \frac{0-1}{1}=-1 \\
				&\left(\epsilon + \frac{\beta\epsilon - \delta\theta}{\la} \right)(1)=(0,\hdots,0)+\frac{(0,\hdots,0)}{1}=(0,\hdots,0) \\
				& \left(\left(\phi + \frac{\beta\phi-\alpha\psi}{\la}\right)(1)\right)^{\text{tr}}=(0,\hdots,0)+\frac{(0,\hdots,0)}{1}=(0,\hdots,0) \\
				&\left(\Xi + \frac{\beta\Xi - \psi\theta}{\lambda}\right)(1)=-I_{(n-2) \times (n-2)} + \frac{0_{(n-2) \times (n-2)}}{1}=-I_{(n-2) \times (n-2)}
			\end{align*}
			
			Similarly, when taking the disjoint union:
			\[\arraycolsep=1pt\def\arraystretch{1.5}
			\begin{array}{c|c}
			\lambda_1 & \; X_1 \; \\
			\hline
			\; X_1 \; & \; A_1 \; \\
			\end{array}
			\cup
			\begin{array}{c|c}
			\lambda_2 & \; X_2 \; \\
			\hline
			\; X_2 \; & \; A_2 \; \\
			\end{array}
			=
			\begin{array}{c|cc}
			\lambda_1 \cdot\lambda_2 & \; X_1 \; & \; X_2 \; \\
			\hline
			X_1 & \lambda_2 A_1 & 0 \\
			X_2 & 0 &  \lambda_1 A_2 \\
			\end{array}\]
			Assuming $\la_1(1)=\la_2(1)=1$ and $A_1(1)=-I_{|X_1| \times |X_1|}, \; A_2(1)=-I_{|X_2| \times |X_2|}$, then for the resulting pair, we have $(\la_1\la_2)(1)=1, \; (\la_2A_1\oplus \la_1A_2)(1)=-I_{(|X_1|+|X_2|)\times (|X_1|+|X_2|)}$.
		\end{proof}
		
		\vspace{0.2cm}
		
		\begin{example} We consider a pure v-tangle $T'$ with diagram $D_{T'}$ given below, which is similar to the one in Example \ref{ex:tangle} (we have omitted the closed component in the earlier tangle $T$ in order to make it a pure tangle). The value of the rMVA invariant is then given as a $3 \times 3$ matrix with an additional element $\la_{T'}$ and a normalizing factor $t^{-1}_1t^{-2}_3$. It can be computed either directly using the degree $0$ and $1$ components of tMVA from the Alexander matrix, or in stages by taking the disjoint union and then gluing together the images of the individual crossings.
			\vspace{0.2cm}
			
			\begin{figure}[H]
				\begin{center}
					\def\svgwidth{8cm}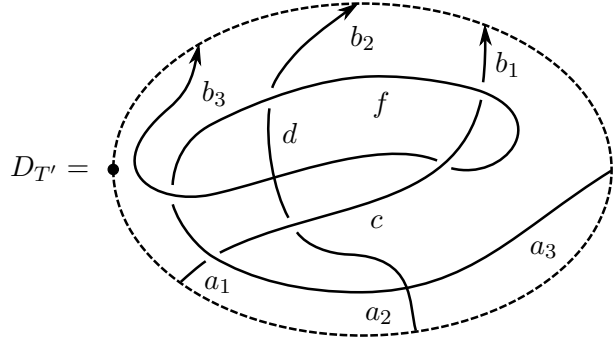
				\end{center}
				\caption{A tangle diagram $D_{T'}$ for the tangle $T'$.}
				\label{fig:rmvaeg1}
			\end{figure}
			
		Then the value of rMVA for the tangle $T'$, $\rMVA(T')$, is given by:
		\begin{figure}[H]
		\begin{center}
		\begin{align*}&\frac{1}{t_1t^2_3}\arraycolsep=1.6pt\def\arraystretch{1.3}
						\begin{array}{c|ccc}
							\la_{T'} &  a_1   &  a_2  &   a_3   \\
							\hline
							b_1 & t_1 t_3^3 \left(t_3 t_1-2 t_1-t_3+1\right)   &   0   & - t_1 t_3^2  \left(t_1-1\right)^2\left(t_3-1\right) \\
							\multirow{2}{*}{$b_2$} & t_3^2\left(t_1 t_3(t_3-2)+\right. & \multirow{2}{*}{$-t_3 \left(t_1+t_3-1\right)$}  & t_3\left(1-t_2\right)\left[t_3+\right.   \\
							& \left.+1-t_3\right)\left(t_2-1\right) & & \left.(t_1-1)(t_1t_3^2-t_1t_3+1)\right] \\
							b_3 & -t_1 t_3^3 \left(t_3-1\right)  & 0  &  t_1 t_3^2 \left(t_3( t_1-2)+1-t_1\right)
						\end{array} \\
						& \\
						&\text{where } \la_{T'}=t_1 t_3^2 \left(t_1+t_3-1\right).\end{align*}
				\end{center}
			\end{figure}
			\label{ex:rangle}
		\end{example}

\section{Metamonoids and $\Gamma$-calculus}
\label{sec:meta}
In this section, we relate the rMVA invariant to a tangle invariant defined by Bar-Natan (\cite{BNS13,DBN1}) which also generalizes the MVA (the generalization is discussed further in Section \ref{sec:ext}). We start by discussing its algebraic structure. Since the gluing operation in the target space of rMVA is not defined on a pair $(\la, \mathcal{A})$ where $\la=0$ (which can occur for tangles with closed components), it does not make it a circuit algebra, and in this setting rMVA is not a circuit algebra morphism. We can however consider rMVA in a different algebraic setting which also describes pure virtual tangles, namely that of {\sl metamonoids}.

\vspace{0.5cm}

\subsection{Metamonoid structure}

\begin{definition}[\cite{BNS13},\cite{DBN1,DBN2,DBN3}]
	A \textbf{metamonoid} is a collection of objects $\{G_X\}$ indexed by finite sets $X$, together with maps between them:
	$$\begin{array}{rll}
	m^{xy}_z&: G_{\{x,y\} \cup X} \longrightarrow G_{\{z\}\cup X} &\text{``multiplication''} \\
	\sqcup&: G_X \times G_Y \longrightarrow G_{X \cup Y} &\text{``union''} \\
	e_x&: G_X \longrightarrow G_{\{x\} \cup X} &\text{``identity''} \\
	\eta_x&: G_{\{x\} \cup X} \longrightarrow G_X &\text{``deletion''} \\
	\sigma^x_z&: G_{\{x\} \cup X} \longrightarrow G_{\{z\}\cup X} &\text{``renaming''}
	\end{array}$$
	satisfying the following relations:

	\vspace{0.2cm}
	
			$\begin{array}{ll}
			&\textbf{Monoid axioms}  \\
			&m^{dc}_{f}\circ m^{ab}_c =m^{cb}_f \circ m^{da}_c \; \text{(associativity)} \\
			&m^{ab}_c\circ e_a =\sigma^b_c  \quad \quad \text{(left identity)} \\
			&m^{ab}_c\circ e_b =\sigma^a_c \quad \quad \text{(right identity)}
			\end{array}$
	
	\vspace{0.2cm}
	
			$\begin{array}{llll}
			&\textbf{Set axioms} & & \\
			&\sigma^a_b \circ e_a =e_b & \sigma^b_c \circ \sigma^a_b = \sigma^a_c & \sigma^b_a \circ \sigma^a_b = \text{Id}    \\
			&\eta_b \circ \sigma^a_b = \eta_a &  \eta_a \circ e_a = \text{Id} &\eta_c \circ m^{ab}_c = \eta_b \circ \eta_a \\
			& \sigma^c_d \circ m^{ab}_c = m^{ab}_d & m^{bc}_d \circ \sigma^a_b = m^{ac}_d &  
			\end{array}$

	\vspace{0.2cm}
	
	\noindent as well as the commuting of operations involving labels which do not interact, for instance $\sigma^a_b \circ \sigma^c_d = \sigma^c_d \circ \sigma^a_b$. Furthermore, the ``union'' operation commutes with all the other operations, for instance we have $\sqcup \circ (m^{xy}_z,m^{uv}_w) = m^{xy}_z \circ m^{uv}_w \circ \sqcup$.
\end{definition}

\begin{remark} An example of a metamonoid, where and element in $G_X$ can be thought of as having a label in $X$, is $G^X=\{f:X \rightarrow G\}$, where $G$ is a group and $X$ is a finite set, with the expected operations. Another example is $pv\mathcal{T}_X$, the collection of $X$-labeled pure regular virtual tangles with the following operations:
	
	$\begin{array}{rll}
	m^{ab}_c&=&\text{join the outgoing end of the strand labelled ``$a$'' to the incoming} \\
	&&\text{end of the strand labelled ``$b$'' and call the resulting strand ``$c$''} \\
	\sqcup &=&\text{take the union of the two tangles by placing them side by side} \\
	e_a&=&\text{add a trivially knotted strand labelled ``$a$''} \\
	\eta_a&=&\text{delete/remove the strand labelled ``$a$'' in the tangle} \\
	\sigma^a_b&=&\text{change the label of strand ``$a$'' to ``$b$''.}
	\end{array}$
	
	Notice in particular that $\eta_x(T) \sqcup \eta_y(T) \neq T \in pv\mathcal{T}_X$, unlike the previous example involving $G^X$.
\end{remark}


The target space of our invariant, together with the gluing and disjoint union maps described in Lemmas \ref{lem:glue} and \ref{lem:union} and the remaining operations defined in a natural way, is also a metamonoid as we discuss below.

\begin{lemma}
	The collection $R_X=\mathbb{R}(\sqrt{t_i}) \times M_{X \times X}(\mathbb{R}(\sqrt{t_i}))$ is a metamonoid under the gluing and disjoint union operations from Section \ref{sec:red}, together with:
	\begin{enumerate}[1)]
		\item
		$\eta_a\left(
		\arraycolsep=6pt\def\arraystretch{1.2}
		\begin{array}{c|cc}
		\la & X & a\\
		\hline
		X & M  & \phi \\
		a & \epsilon & \alpha
		\end{array}
		\right) = \left(\arraycolsep=6pt\def\arraystretch{1.2}\begin{array}{c|c}
		\la & X \\
		\hline
		X & M
		\end{array}\right)_{t_a=1}$
		
		\vspace{0.2cm}
		
		\item $\sigma^a_b\left(
		\arraycolsep=6pt\def\arraystretch{1.2}
		\begin{array}{c|cc}
		\la & X & a\\
		\hline
		X & M  & \phi \\
		a & \epsilon & \alpha
		\end{array}
		\right) =\left(\arraycolsep=6pt\def\arraystretch{1.2}\begin{array}{c|cc}
		\la & X & b\\
		\hline
		X & M  & \phi \\
		b & \epsilon & \alpha
		\end{array}\right)_{t_a\rightarrow t_b}$
	
		\item $e_a\left(
		\arraycolsep=4pt\def\arraystretch{1.2}
		\begin{array}{c|c}
		\la & X \\
		\hline
		X & M
		\end{array}\right) =
		\arraycolsep=6pt\def\arraystretch{1.2}
		\begin{array}{c|cc}
		\la & X & a\\
		\hline
		X & M  & 0 \\
		a & 0 & -\la
		\end{array}$
	\end{enumerate}
	\label{lem:rxmet}
\end{lemma}

\begin{proof}
	The commuting and set axioms follow directly from the definition of the operations on $R_X$, so we only need to verify that the monoid axioms are satisfied:
	
	\vspace{0.5cm}
	
	\underline{Checking the ``left identity'' axiom:}
	\def\arraystretch{1.2}
	\begin{align*}
		&m^{ab}_c \circ e_a \left(
		\begin{array}{c|cc}
			\la & b & X \\
			\hline
			b & \delta & \epsilon \\
			X & \psi & \Xi
		\end{array} \right) =
		m^{ab}_c\left(\begin{array}{c|ccc}
			\la & a & b & X \\
			\hline
			a & -\la & 0 & 0_{|X|} \\
			b & 0 & \delta & \epsilon \\
			X & 0^{\text{tr}}_{|X|} & \psi & \Xi
		\end{array}
		\right) \\
		&=\left({\setstretch{1.5}\begin{array}{c|cc}
				\la & c & X \\
				\hline
				c & 0+\frac{0-\delta(-\la)}{\la} & \epsilon+\frac{0\cdot \epsilon-\delta\cdot 0_{|X|}}{\la}  \\
				X & 0^{\text{tr}}_{|X|} +\frac{0^{\text{tr}}_{|X|}-(-\la)\psi}{\la} & \Xi+\frac{0\cdot \Xi - \psi \cdot 0_{|X|}}{\la}
			\end{array}}\right)_{\hspace{-0.15cm}t_a,t_b\rightarrow t_c} \\
			&=\left(\begin{array}{c|cc}
				\la & c & X \\
				\hline
				c & \delta & \epsilon \\
				X & \psi & \Xi
			\end{array}\right)_{\hspace{-0.15cm}t_b\rightarrow t_c}\hspace{-0.7cm}=\sigma^b_c\left(\begin{array}{c|cc}
			\la & b & X \\
			\hline
			b & \delta & \epsilon \\
			X & \psi & \Xi
		\end{array}\right)
	\end{align*}
	
	\vspace{0.5cm}
	
	\underline{Checking the ``right identity'' axiom:}
	\def\arraystretch{1.2}
	
		\begin{align*}&m^{ab}_c \circ e_b \left(
		\begin{array}{c|cc}
			\la & a & X \\
			\hline
			a & \alpha & \theta \\
			X & \phi & \Xi
		\end{array} \right) =
		m^{ab}_c\left(\begin{array}{c|ccc}
			\la & a & b & X \\
			\hline
			a & \alpha & 0 & \theta \\
			b & 0 & -\la &  0_{|X|} \\
			X & \phi & 0^{\text{tr}}_{|X|} & \Xi
		\end{array}
		\right)
		\end{align*}
		
		\begin{align*}
		&=\left({\setstretch{1.5}\begin{array}{c|cc}
				\la & c & X \\
				\hline
				c & 0+\frac{0-\alpha(-\la)}{\la} & 0_{|X|}+\frac{0_{|X|}-(-\la)\theta}{\la}  \\
				X & \phi +\frac{0 \cdot \phi -\alpha \cdot 0^{\text{tr}}_{|X|}}{\la} & \Xi+\frac{0\cdot \Xi - 0^{\text{tr}}_{|X|} \cdot \theta}{\la}
			\end{array}}\right)_{\hspace{-0.15cm}t_a,t_b\rightarrow t_c} \\
		&=\left(\begin{array}{c|cc}
				\la & c & X \\
				\hline
				c & \alpha & \theta \\
				X & \phi & \Xi
			\end{array}\right)_{\hspace{-0.15cm}t_a\rightarrow t_c}\hspace{-0.7cm}=\sigma^a_c\left(\begin{array}{c|cc}
			\la & a & X \\
			\hline
			a & \alpha & \theta \\
			X & \phi & \Xi
		\end{array}\right)
	\end{align*}
	
	\vspace{0.2cm}
	
	\underline{Checking the ``associativity'' axiom:}
	\def\arraystretch{1.2}
	\begin{align*}
		&m^{dc}_f \circ m^{ab}_c\left(\begin{array}{c|cccc}
			\la & a & b & d & X \\
			\hline
			a & \alpha & \beta & \tau & \theta \\
			b & \gamma & \delta & \nu & \ep \\
			d & \kappa & \rho & \xi & \sigma \\
			X & \phi & \psi & \mu & \Xi
		\end{array}\right) = \\
		&=m^{dc}_f\left({\setstretch{1.5}
			\begin{array}{c|ccc}
				\la +\beta & c & d & X \\
				\hline
				c & \gamma+\frac{\beta\gamma-\alpha\delta}{\la} & \nu + \frac{\beta\nu - \delta\tau}{\la} & \ep + \frac{\beta\ep - \delta\theta}{\la} \\
				d & \kappa + \frac{\beta\kappa-\alpha\rho}{\la} & \xi + \frac{\beta\xi-\rho\tau}{\la} & \sigma + \frac{\beta\sigma - \rho\theta}{\la} \\
				X & \phi + \frac{\beta\phi - \alpha\psi}{\la} & \mu + \frac{\beta\mu-\tau\psi}{\la} & \Xi + \frac{\beta\Xi-\psi\theta}{\la}
			\end{array}} \right)_{t_a,t_b \rightarrow t_c} \\
			& \\
			&=\left({\setstretch{1.5}
				\begin{array}{c|cc}
					\frac{(\beta +\la) (\kappa +\la)-\alpha  \rho }{\la} & f & X \\
					\hline
					f & \mathcal{A}_{ff} & \mathcal{A}_{fX} \\
					X & \mathcal{A}_{Xf} & \mathcal{A}_{XX}
				\end{array}} \right)_{t_c,t_d \rightarrow t_f}
			\end{align*}
			Where:
			\begin{align*}
				\mathcal{A}_{ff} &=\frac{\alpha(  \delta  \xi - \rho\nu) -(\beta + \la)\gamma  \xi + (\kappa +\lambda )((\beta + \la)\nu-\tau\delta )  + \rho  \tau\gamma}{\la^2} \\
				\mathcal{A}_{fX} &=\frac{\alpha  (\delta  \sigma -\rho  \epsilon )-(\beta +\la)\gamma  \sigma  +(\kappa +\lambda ) ((\beta +\lambda )\epsilon  -\delta  \theta )+\rho\gamma\theta  }{\la^2} \\
				\mathcal{A}_{Xf} &=\frac{\alpha  (\xi  \psi -\rho \mu )-(\beta +\lambda )\xi  \phi  +(\kappa +\lambda ) ( (\beta +\lambda )\mu -\tau  \psi )+\rho  \tau  \phi}{\la^2} \\
				\mathcal{A}_{XX} &=\frac{\alpha  (\psi \sigma -  \rho \Xi)+(\kappa +\lambda ) ((\beta +\lambda )\Xi  - \psi \theta )-(\beta +\lambda )\phi \sigma  +\rho  \phi\theta  }{\la^2}
			\end{align*}
			
			Similarly, for the right side of the axiom we get:
			
			\def\arraystretch{1.2}
			\begin{align*}
				&m^{cb}_f \circ m^{da}_c\left(\begin{array}{c|cccc}
					\la & a & b & d & X \\
					\hline
					a & \alpha & \beta & \tau & \theta \\
					b & \gamma & \delta & \nu & \ep \\
					d & \kappa & \rho & \xi & \sigma \\
					X & \phi & \psi & \mu & \Xi
				\end{array}\right) = \\
				&=m^{cb}_f\left({\setstretch{1.5}
					\begin{array}{c|ccc}
						\la +\kappa & c & b & X \\
						\hline
						c & \tau+\frac{\kappa\tau-\alpha\xi}{\la} & \beta + \frac{\kappa\beta-\alpha\rho}{\la} & \theta + \frac{\kappa\theta-\alpha\sigma}{\la} \\
						b & \nu + \frac{\kappa\nu-\xi\gamma}{\la} & \delta + \frac{\kappa\delta-\gamma\rho}{\la} & \ep + \frac{\kappa\ep-\gamma\sigma}{\la} \\
						X & \mu + \frac{\kappa\mu-\xi\phi}{\la} & \psi + \frac{\kappa\psi-\rho\phi}{\la} & \Xi + \frac{\kappa\Xi-\phi\sigma}{\la}
					\end{array}} \right)_{t_a,t_d \rightarrow t_c} \\
					&=\left({\setstretch{1.5}
						\begin{array}{c|cc}
							\frac{(\beta +\la) (\kappa +\la)-\alpha  \rho }{\la} & f & X \\
							\hline
							f & \mathcal{B}_{ff}  & \mathcal{B}_{fX} \\
							X & \mathcal{B}_{Xf} & \mathcal{B}_{XX}
						\end{array}} \right)_{t_c,t_b \rightarrow t_f}
					\end{align*}
					Where:
					\begin{align*}
					\mathcal{B}_{ff} &=\frac{\alpha(  \delta  \xi - \rho\nu) -(\beta + \la)\gamma  \xi + (\kappa +\lambda )((\beta + \la)\nu-\tau\delta )  + \rho  \tau\gamma}{\la^2} \\
					\mathcal{B}_{fX} &=\frac{\alpha  (\delta  \sigma -\rho  \epsilon )-(\beta +\la)\gamma  \sigma  +(\kappa +\lambda ) ((\beta +\lambda )\epsilon  -\delta  \theta )+\rho\gamma\theta  }{\la^2} \\
					\mathcal{B}_{Xf} &=\frac{\alpha  (\xi  \psi -\rho \mu )-(\beta +\lambda )\xi  \phi  +(\kappa +\lambda ) ( (\beta +\lambda )\mu -\tau  \psi )+\rho  \tau  \phi}{\la^2} \\
					\mathcal{B}_{XX} &=\frac{\alpha  (\psi \sigma -  \rho \Xi)+(\kappa +\lambda ) ((\beta +\lambda )\Xi  - \psi \theta )-(\beta +\lambda )\phi \sigma  +\rho  \phi\theta  }{\la^2}
					\end{align*}
				\end{proof}
				
				Thus, both the domain $pv\mathcal{T}_X$ and the target space $R_X$ of rMVA are metamonoids. Furthermore, in this more general algebraic setting, the tangle invariant $\rMVA: pv\mathcal{T}_X \longrightarrow R_X$ is a \textbf{metamonoid morphism} by construction. We will next relate it to a tangle invariant defined by Bar-Natan in this algebraic context.
				
				\subsection{$\Gamma$- or Gassner-calculus}
				
				In \cite{DBN1} and \cite{DBN2}, Bar-Natan defines another tangle invariant, ``$\Gamma$-calculus'', which also generalizes the multivariable Alexander polynomial. It is also an invariant for pure virtual regular $X$-labelled tangles, where $X$ is a finite set. Hence, its domain is the metamonoid $pv\mathcal{T}_X$. The target space is the metamonoid $\Gamma_X=\mathbb{Z}(t_i)_{i \in X} \times M_{X \times X}(\mathbb{Z}(t_i))$ with operations given by:
				
				\vspace{0.5cm}
				
				\textbf{Multiplication:}
				\[\arraycolsep=5pt
				\begin{array}{c|ccc}
				\lambda & a & b & X \\
				\hline
				a & \alpha & \beta & \theta \\
				b & \gamma & \delta & \epsilon \\
				X & \phi & \psi & \Xi 
				\end{array}
				\xrightarrow{\; m^{a,b}_c \;}
				\left(\arraycolsep=10pt\def\arraystretch{1.5}
				\begin{array}{c|ccc}
				\lambda(1-\beta) & c & X \\
				\hline
				c & \gamma + \frac{\alpha\delta}{1-\beta} & \epsilon + \frac{\delta\theta}{1-\beta} \\
				X & \phi + \frac{\alpha\psi}{1-\beta} & \Xi + \frac{\psi\theta}{1-\beta} 
				\end{array}\right)_{t_a,t_b \rightarrow t_c}\]
				\vspace{0.2cm}
				
				\begin{remark}Analogously to the case of $R_X$, multiplication in $\Gamma_X$ is well-defined on the image of pure tangles. Indeed, for a pure tangle $T$ whose image in $\Gamma_X$ is the pair $(\la_T,\mathcal{A}_T)$, we have $\la_T(1)=1, \mathcal{A}_T(1)=\text{Id}$. Namely, setting all the variables $t_i=1$, we get the pair $(1,\text{Id})$. So, in the multiplication operation, the division term $1-\beta$ where $\beta$ is an off-diagonal element is never identically zero.
				\end{remark}
				\vspace{0.2cm}
				
				\begin{flushleft}
					\begin{minipage}[t]{0.58\textwidth}
						\textbf{Union:}
						\[\arraycolsep=1pt\def\arraystretch{1.2}
						\begin{array}{c|c}
						\lambda_1 & \; X_1 \; \\
						\hline
						\; X_1 \; & \; A_1 \; \\
						\end{array}
						\sqcup
						\begin{array}{c|c}
						\lambda_2 & \; X_2 \; \\
						\hline
						\; X_2 \; & \; A_2 \; \\
						\end{array}
						=\arraycolsep=2pt\def\arraystretch{1.2}
						\begin{array}{c|cc}
						\lambda_1 \cdot\lambda_2 & \; X_1 \; & \; X_2 \; \\
						\hline
						X_1 & A_1  & 0 \\
						X_2 & 0 &  A_2 \\
						\end{array}\]
					\end{minipage}\hspace{0.5cm}
					\begin{minipage}[t]{0.3\textwidth}
						\textbf{Identity:}
						\[
						e_a\left(
						\arraycolsep=2pt\def\arraystretch{1.2}
						\begin{array}{c|c}
						\la & X \\
						\hline
						X & A
						\end{array}\right)=
						\arraycolsep=4pt\def\arraystretch{1.2}\begin{array}{c|cc}
						\la & X & a\\
						\hline
						X & A  & 0 \\
						a & 0 & 1
						\end{array}\]
					\end{minipage}
				\end{flushleft}
				\vspace{0.2cm}
				
				The deletion and renaming operations for $\Gamma_X$ are the same as those for $R_X$ and we will not repeat them here. The Bar-Natan invariant is then defined as the metamonoid morphism:
				\[Z: pv\mathcal{T}_X \longrightarrow \Gamma_X\]
				which maps the positive and negative crossing generators as indicated below.
					
					\begin{figure}[H]
					\begin{center}
						\begin{minipage}[c]{0.12\textwidth}
							\def\svgwidth{1.4cm}\input{rPosX.pdf_tex}
						\end{minipage}
						\begin{minipage}[c]{0.3\textwidth}
							$\xrightarrow{\quad Z \quad} 
							\arraycolsep=4pt\def\arraystretch{1.3}
							\begin{array}{c|cc}
							1 & a & b \\
							\hline
							a & 1 &  1-t_a \\
							b & 0 & t_a 
							\end{array}$
						\end{minipage}
						\hspace{1cm}
						\begin{minipage}[c]{0.12\textwidth}
							\def\svgwidth{1.4cm}\input{rNegX.pdf_tex}
						\end{minipage}
						\begin{minipage}[c]{0.3\textwidth}
							$\xrightarrow{\quad Z \quad}
							\arraycolsep=4pt\def\arraystretch{1.3} 
							\begin{array}{c|cc}
							1 & a & b \\
							\hline
							a & 1 &  1-t^{-1}_a \\
							b & 0 & t^{-1}_a 
							\end{array}$
						\end{minipage}
					\end{center}
					\caption{The images under $Z$ of $\overcrossing$ and $\undercrossing$.}
					\end{figure}

				\subsection{Relating rMVA and $Z$}
				Both tangle invariants rMVA and $Z$ are metamonoid morphisms mapping into almost identical spaces, so unsurprisingly there is a close relationship between them. We will explore that relationship next by first defining a slight variation of Bar-Natan's invariant. Note that the elements in both $R_X$ and $\Gamma_X$ are matrices of size $|X| \times |X|$ with an additional element, and if we equip $R_X$ with the operations from $\Gamma_X$ instead, it is still a metamonoid. So, we get a second metamonoid structure on $R_X$, which we will denote by $\widetilde{R}_X$. We consider a slight modification $\widetilde{Z}$  of Bar-Natan's $Z$ invariant, namely the metamonoid morphism:
				\[\widetilde{Z}: pv\mathcal{T}_X \longrightarrow  \widetilde{R}_X\]
				
				\noindent which is defined on the generators as:
				
				\begin{figure}[H]
					\begin{center}
						\begin{minipage}[c]{0.12\textwidth}
							\def\svgwidth{1.4cm}\input{rPosX.pdf_tex}
						\end{minipage}
						\begin{minipage}[c]{0.3\textwidth}
							$\xrightarrow{\quad \widetilde{Z} \quad}
							\arraycolsep=4pt\def\arraystretch{1.3}
							\begin{array}{c|cc}
							t^{\frac{1}{2}}_a & a & b \\
							\hline
							a & 1 &  0 \\
							b & \frac{t_b-1}{t_a} & t^{-1}_a 
							\end{array}$
						\end{minipage}
						\hspace{1cm}
						\begin{minipage}[c]{0.12\textwidth}
							\def\svgwidth{1.4cm}\input{rNegX.pdf_tex}
						\end{minipage}
						\begin{minipage}[c]{0.32\textwidth}
							$\xrightarrow{\quad \widetilde{Z} \quad}
							\arraycolsep=4pt\def\arraystretch{1.3} 
							\begin{array}{c|cc}
							t^{-\frac{1}{2}}_a & a & b \\
							\hline
							a & 1 &  0 \\
							b & 1-t_b & t_a 
							\end{array}$
						\end{minipage}
					\end{center}
					\caption{The images under $\widetilde{Z}$ of $\overcrossing$ and $\undercrossing$.}
					\label{fig:wzgen}
				\end{figure}
				Note that, analogously to the previous cases, for each of these pairs $(\la,\mathcal{A})$, setting all the variables $t_i$ to $1$ gives $(\la(1),\mathcal{A}(1))=(1,\text{Id})$. So, for any pure tangle $T$, we have $\la_T \neq 0$.
				This version of the Bar-Natan invariant turns out to be equivalent to rMVA, as we show next. 
				
				\begin{theorem}
					There is a metamonoid morphism $\widetilde{R}_X \rightarrow R_X$, which is an isomorphism on the level of the images of rMVA and $\widetilde{Z}$, taking the positive and negative crossing generators for $\widetilde{Z}$ in $\widetilde{R}_X$ to those for rMVA in $R_X$.
					\label{thm:iso}
				\end{theorem}

				\begin{proof}
					The map $F: \widetilde{R}_X \longrightarrow R_X$ defined by $F(\la,\mathcal{A})=(\la,-\la \cdot \mathcal{A})$ takes the positive and negative crossing generators for $\widetilde{Z}$ in Figure \ref{fig:wzgen} to precisely the positive and negative crossing generators for rMVA in Figure \ref{fig:pnmat}. Furthermore, for any pair $(\la,\mathcal{A})$ in the image $\widetilde{Z}(pv\mathcal{T}_X)$ or $\rMVA(pv\mathcal{T}_X)$ of pure tangles, we have seen that $\la \neq 0$. So, the restriction of $F$ to these subsets is a bijection. It remains to show that $F$ is a metamonoid morphism, i.e. we need to check that the following diagram commutes, as well as the analogous ones with $m^{ab}_c$ replaced by $\sqcup, \; e_a, \; \eta_a, \; \sigma^a_b$.
					\[
					\begin{tikzcd}[ampersand replacement=\&]
					\arraycolsep=5pt
					\widetilde{R}_{X\cup\{a,b\}} 
					\arrow{r}{F} \arrow[swap]{d}{(m^{ab}_c)_{\widetilde{R}_X}} \& 
					R_{X\cup\{a,b\}}  \arrow{d}{(m^{ab}_c)_{R_X}} \\
					\widetilde{R}_{X\cup\{c\}}  \arrow{r}{F} \& 
					R_{X\cup\{c\}} 
					\end{tikzcd}
					\]
					
					For a given element in $\widetilde{R}_{X\cup\{a,b\}}$, we first follow the diagram horizontally then vertically:
					
					\begin{align*}
						&(m^{ab}_c)_{R_X}\circ F\left(\arraycolsep=5pt
						\begin{array}{c|ccc}
							\la & a & b & X \\
							\hline
							a & \alpha & \beta & \theta \\
							b & \gamma & \delta & \epsilon \\
							X & \phi & \psi & \Xi 
						\end{array} \right) = (m^{ab}_c)_{R_X} \left(
						\begin{array}{c|ccc}
							\la & a & b & X \\
							\hline
							a & -\la\alpha & -\la\beta & -\la\theta \\
							b & -\la\gamma & -\la\delta & -\la\epsilon \\
							X & -\la\phi & -\la\psi & -\la\Xi 
						\end{array} \right)
					\end{align*}
						\begin{align*}
						&=\left(\arraycolsep=4pt\def\arraystretch{1.5}
						\begin{array}{c|cc}
							\la(1-\beta) & c & X \\
							\hline
							c & -\la\gamma + \la\beta\gamma-\la\delta\alpha & -\la\ep + \la\beta\ep-\la\delta\theta \\
							X & -\la\phi + \la\beta\phi-\la\alpha\psi & -\la\Xi + \la\beta\Xi-\la\psi\theta
						\end{array} \right)_{t_a,t_b \rightarrow t_c} \\
						& \\
						&=\left(\arraycolsep=4pt\def\arraystretch{1.5}
						\begin{array}{c|ccc}
							\lambda(1-\beta) & c & X \\
							\hline
							c & -\lambda(1-\beta)\gamma -\la\alpha\delta & -\lambda(1-\beta)\epsilon -\la\delta\theta \\
							X & -\lambda(1-\beta)\phi  -\la\alpha\psi & -\lambda(1-\beta)\Xi -\la\psi\theta 
						\end{array}\right)_{t_a,t_b \rightarrow t_c} 
					\end{align*}
					
					Alternatively, following the map going down first, then right, we get:
					
					\begin{align*}
						&F \circ(m^{ab}_c)_{\widetilde{R}_X}\left(\arraycolsep=3pt
						\begin{array}{c|ccc}
							\la & a & b & X \\
							\hline
							a & \alpha & \beta & \theta \\
							b & \gamma & \delta & \epsilon \\
							X & \phi & \psi & \Xi 
						\end{array} \right) =F\left(\arraycolsep=3pt\def\arraystretch{1.2}
						\begin{array}{c|ccc}
							\lambda(1-\beta) & c & X \\
							\hline
							c & \gamma + \frac{\alpha\delta}{1-\beta} & \epsilon + \frac{\delta\theta}{1-\beta} \\
							X & \phi + \frac{\alpha\psi}{1-\beta} & \Xi + \frac{\psi\theta}{1-\beta} 
						\end{array}\right)_{t_a,t_b \rightarrow t_c} \\
						& \\
						&=\left(\arraycolsep=4pt\def\arraystretch{1.5}
						\begin{array}{c|ccc}
							\lambda(1-\beta) & c & X \\
							\hline
							c & -\lambda(1-\beta)\gamma -\la\alpha\delta & -\lambda(1-\beta)\epsilon -\la\delta\theta \\
							X & -\lambda(1-\beta)\phi  -\la\alpha\psi & -\lambda(1-\beta)\Xi -\la\psi\theta 
						\end{array}\right)_{t_a,t_b \rightarrow t_c} 
					\end{align*}
					
					The proof that the diagram commutes for $\eta_a$ and $\sigma^a_b$ is immediate from the definition, so we will additionally only check it for $\sqcup$ and $e_a$, whose verification is also straightforward. The following diagram verifies that $F$ preserves disjoint unions:
					
					\[
					\begin{tikzcd}[ampersand replacement=\&]
					\arraycolsep=3pt\def\arraystretch{1.2}
					\left(\begin{array}{c|c}
					\la_1 & \; X_1 \; \\
					\hline
					\; X_1 \; & \; A_1 \; \\
					\end{array}\;,\;
					\begin{array}{c|c}
					\la_2 & \; X_2 \; \\
					\hline
					\; X_2 \; & \; A_2 \; \\
					\end{array} \right)
					\arrow{r}{F} \arrow[swap]{d}{(\sqcup)_{\widetilde{R}_X}} \& 
					\arraycolsep=3pt\def\arraystretch{1.2}
					\left(\begin{array}{c|c}
					\la_1 & \; X_1 \; \\
					\hline
					\; X_1 \; & \; -\la_1 A_1 \; \\
					\end{array}\;,\;
					\begin{array}{c|c}
					\lambda_2 & \; X_2 \; \\
					\hline
					\; X_2 \; & \; -\la_2 A_2 \; \\
					\end{array} \right)  \arrow{d}{(\sqcup)_{R_X}} \\
					\arraycolsep=3pt\def\arraystretch{1.2}
					\begin{array}{c|cc}
					\la_1 \cdot\la_2 & \; X_1 \; & \; X_2 \; \\
					\hline
					X_1 & A_1  & 0 \\
					X_2 & 0 &  A_2 \\
					\end{array}  \arrow{r}{F} \& 
					\arraycolsep=3pt\def\arraystretch{1.2}
					\begin{array}{c|cc}
					\la_1 \cdot\la_2 & \; X_1 \; & \; X_2 \; \\
					\hline
					X_1 & - \la_1\la_2 A_1  & 0 \\
					X_2 & 0 &  -\la_1\la_2 A_2 \\
					\end{array} 
					\end{tikzcd}
					\]
					Similarly, we confirm with the diagram below that $F$ commutes with the ``identity'' map $e_a$ of the two metamonoids:
					
					\[
					\begin{tikzcd}[ampersand replacement=\&]
					\arraycolsep=3pt\def\arraystretch{1.2}
					\begin{array}{c|c}
					\la & \; X \; \\
					\hline
					\; X \; & \; A \; \\
					\end{array}
					\arrow{r}{F} \arrow[swap]{d}{(e_a)_{\widetilde{R}_X}} \& 
					\arraycolsep=3pt\def\arraystretch{1.2}
					\begin{array}{c|c}
					\la & \; X \; \\
					\hline
					\; X \; & \; -\la A \; \\
					\end{array} \arrow{d}{(e_a)_{R_X}} \\
					\arraycolsep=3pt\def\arraystretch{1.2}
					\begin{array}{c|cc}
					\la & \; X \; & \; a \; \\
					\hline
					X & A  & 0 \\
					a & 0 &  1 \\
					\end{array}  \arrow{r}{F} \& 
					\arraycolsep=3pt\def\arraystretch{1.2}
					\begin{array}{c|cc}
					\la & \; X \; & \; a \; \\
					\hline
					X & - \la A  & 0 \\
					a & 0 &  -\la \\
					\end{array} 
					\end{tikzcd}
					\]
				\end{proof}

					\begin{remark}
						Let $P_X=\rMVA(pv\mathcal{T}_X)$ and $\widetilde{P}_X=\widetilde{Z}(pv\mathcal{T}_X)$ denote the image sets of the two invariants. Then we can summarize the main result from Theorem \ref{thm:iso} in Figure \ref{fig:mmcom}, namely that the target spaces of rMVA and $\widetilde{Z}$ are isomorphic and the isomorphism takes the images of the tangle metamonoid generators using $\widetilde{Z}$ to those using rMVA. In particular, this tells us that $\rMVA$ and $\widetilde{Z}$ are equivalent as pure tangle invariants.
					\end{remark}
		
						\begin{figure}[H]
							\begin{center}$\xymatrix{
									& pv\mathcal{T}_X \ar[ld]_{\widetilde{Z}} \ar[rd]^{\rMVA} &\\
									\widetilde{R}_X \supseteq \widetilde{P}_X \ar[rr]^{\cong}_{F} & & P_X \subseteq R_X }$
							\end{center}
							\captionsetup{width=0.9\textwidth}
							\caption{A commutative diagram of metamonoid morphisms.}
							\label{fig:mmcom}
						\end{figure}

	\section{Links, Braids, and Tangles with Closed Components}
\label{sec:braids}

\subsection{Recovering the MVA}
\label{sec:MVA} One of the main motivations for studying the rMVA tangle invariant is that, since it comes from Archibald's tMVA invariant, it is a generalization to tangles of the multivariable Alexander polynomial on links defined by Torres, \cite{Tor53}. It is therefore of interest to study how to recover the MVA directly from the rMVA.  

In order to consider links, we need to look at tangles with closed components, so setting aside the algebraic structure of metamonoids, let us consider the larger space $v\mathcal{T}_X$ of (not necessarily pure) regular virtual tangles labelled using the set $X$. Note that we can also define rMVA on $v\mathcal{T}_X$, not in terms of a metamonoid morphism presented in Lemmas \ref{lem:glue}, \ref{lem:union}, and \ref{lem:rxmet}, but directly from the degree $0$ and $1$ components using the Hodge star operator on the image of tMVA (which can also be found through the Alexander matrix) and the normalizing factor, as in Equation \ref{eqn:dirdef}. More precisely, for a tangle $T \in v\mathcal{T}_X$ possibly having some closed components, we can define: \[\rMVA(T):=\prod_k{t^{-\frac{\mu(k)}{2}}_k}(\la_T,\mathcal{A}_T)\vspace{-0.15cm}\] In this setting we no longer have the gluing operation. This has the disadvantage of no longer allowing the ``divide-and-conquer'' approach. Namely, computing the invariant by splitting a tangle into pieces through the generators and then gluing in the context of metamonoid morphisms. However, it allows us to consider the invariant on a bigger space which includes tangles with closed components, and to recover from it the multivariable Alexander polynomial for links.

The MVA on links is in fact defined using the same construction for the Alexander matrix as described in Figure \ref{fig:rules}. The version of the link invariant that we will work with and for which it is simplest to describe the relation to rMVA is \textbf{vMVA}, or the multivariable Alexander polynomial for regular long virtual links, i.e. links in which exactly one strand has both ends going to infinity.

\begin{definition}[\cite{Arch10}]
	The multivariable Alexander polynomial vMVA for a regular long virtual link $L$ is defined as:
	\[\vMVA(L)=\frac{1}{t_{l}-1}\prod_s{t^{-\frac{\mu(s)}{2}}_s}\det{M(D_L)^{l^{\text{out}};\widehat{l^{\text{in}}}}}\]
	where $t_l$ is the variable associated to the long strand, $\mu(s)$ as before is the number of times the strand labelled ``$s$'' overcrosses in a crossing, $D_L$ is a diagram of the link $L$, and $M(D_L)^{l^{\text{out}};\widehat{l^{\text{in}}}}$ is the submatrix of the Alexander matrix with all the columns corresponding to the internal arcs included, as well as the column corresponding to the outgoing arc of strand ``$l$'', but not the incoming one.
\end{definition}

Going back to our consideration of rMVA, another disadvantage of using it on the bigger space $v\mathcal{T}_X$ as above is that for a tangle $T \in v\mathcal{T}_X$, we are no longer guaranteed that $\la_T \neq 0$, and so we cannot apply Theorem \ref{thm:deg01} to conclude that $\rMVA(T)$ determines $\tMVA(T)$. We can however say more if we restrict our attention to the space $v\mathcal{T}^1_X$ of $1$-tangles, i.e. tangles with possibly some closed components and exactly one open strand component, where the set $X$ labels both the closed components and the open component. In that case, the Alexander matrix has a single outgoing arc column and a single incoming arc column, so $\tMVA(v\mathcal{T}^1_X)$ lives in only the degree $0$ and $1$ components and is therefore completely determined by rMVA directly from the definition. For $T \in v\mathcal{T}^1_X$:
\[\rMVA\left(\xrightarrow{a^{\text{in}}}\hspace{-0.14cm}\fbox{T}\hspace{-0.14cm}\xrightarrow{a^{\text{out}}}\right)=
\arraycolsep=3pt\def\arraystretch{1.2}
\begin{array}{c|c} 
\la_T & a^{\text{in}} \\ 
\hline
a^{\text{out}} & \mathcal{A}_T
\end{array}\]
where the degree zero element is $\la_T=\prod_k{t^{-\frac{\mu(k)}{2}}_k}\det{M(D_T)^{a^{\text{out}};\emptyset}}$, and the matrix of degree one components is given by  $\mathcal{A}_T=\prod_k{t^{-\frac{\mu(k)}{2}}_k}\det{M(D_T)^{\emptyset;a^{\text{in}}}}$. Thus, 
\begin{align*}
	\tMVA(T)&=\prod_k{t^{-\frac{\mu(k)}{2}}_k} a^{\text{out}} \otimes \left( \det{M(D_T)^{a^{\text{out}};\emptyset}}a^{\text{out}} + \det{M(D_T)^{\emptyset;a^{\text{in}}}}a^{\text{in}}\right) \\
	&=a^{\text{out}}\otimes \left(\la_T \; a^{\text{out}} + \mathcal{A}_T \; a^{\text{in}} \right).
\end{align*}

To make the connection to links, we note that long links, i.e. those with exactly one component having ends going to infinity, coincide with $1$-tangles. So, we can relate the values of $\rMVA$ on $1$-tangles and $\vMVA$ on long links. Archibald makes this connection between $\tMVA$ and $\vMVA$ in \cite{Arch10} by noting that the Alexander matrix for the two invariants is constructed by the same rules. Furthermore, as can be deduced from Figure \ref{fig:rules}, the columns $\{C_s\}$ of the matrix satisfy the relation $\sum_s(t_s-1)C_s=0$. Since the incoming and outgoing ends of a tangle correspond to the same strand and so to the same variable $t_s$, the relation is:
\begin{theorem}{\cite[Theorem 6.8]{Arch10}}
	For a $1$-tangle $T$: \[\frac{1}{t_a-1}\tMVA\left(\xrightarrow{a^{\text{in}}}\hspace{-0.14cm}\fbox{T}\hspace{-0.14cm}\xrightarrow{a^{\text{out}}} \right)=\vMVA\left(\xrightarrow{a^{\text{in}}}\hspace{-0.14cm}\fbox{T}\hspace{-0.14cm}\xrightarrow{a^{\text{out}}} \right)a^{\text{out}}\otimes ( a^{\text{out}}- a^{\text{in}}) \]
\end{theorem}
This theorem relies in particular on a conclusion based on the relation among the columns of the Alexander matrix, namely that $\det{M(D_T)^{a^{\text{out}};\emptyset}}=-\det{M(D_T)^{\emptyset;a^{\text{in}}}}$. This leads to the connection between $\rMVA$ and $\vMVA$:

\begin{proposition}
	For $\xrightarrow{a^{\text{in}}}\hspace{-0.14cm}\fbox{T}\hspace{-0.14cm}\xrightarrow{a^{\text{out}}}$ a regular virtual $1$-tangle (i.e. a long knot or link) and $\rMVA(T)=(\la_T,\mathcal{A}_T)$ (with the normalizing factor implicitly included):
	\[\la_T=-\mathcal{A}_T \quad \quad \frac{\la_T}{t_a-1}=\vMVA(T)\]
\end{proposition}

For long knots, i.e. $1$-tangles without closed components, we can still benefit from the metamonoid structure in using this proposition to break down the knot into smaller pieces and then glue their images, making computations more efficient.

\subsection{Braids and the Gassner representation}

Another specialization of the rMVA invariant leads again to familiar terrirory. For this section we will restrict our attention to a subset of the collection of $X$-labelled tangles, namely that of pure virtual braids. Note that ``pure'' in this case refers to the more standard meaning of the underlying permutation being the identity.

\begin{definition}[\cite{Ba04}\cite{BND14}]
	The pure virtual braid group $pv\mathcal{B}_X$ with strands labelled by $X$ is generated by $\si_{ab}$ for all $a\neq b \in X$, corresponding to the positive crossing where strand ``$a$'' crosses over strand ``$b$'', and relations:
	\begin{align*}\si_{ab}\si_{ac}\si_{bc} &=\si_{bc}\si_{ac}\si_{ab} \quad \quad \forall \; a,b,c \text{ pairwise distinct} \\
		\si_{ab}\si_{cd}&=\si_{cd}\si_{ab} \quad \quad \forall \; a,b,c,d \text{ pairwise distinct}
	\end{align*}
\end{definition}
\begin{remark}
	Note that in this notation, the labels record the ``identity'' of the strand rather than its location. This is unlike the usual braid group where only adjacent strands can cross and two strands crossing exchange their labels (which indicate their position). In the notation above for this pure virtual setting, a strand labelled ``$a$'' continues to carry that label even after crossing with another strand.
\end{remark}
In her thesis \cite{G59}, Gassner defines a multivariable version of the Burau representation for the braid group, which can be generalized to $pv\mathcal{B}_X$ by mapping each generator $\si_{ab}, \; a,b \in X$, to the following matrix (See also \cite{DBN4}):
\[\si_{ab} \xrightarrow{G} \left[
\begin{array}{c|ccc}
& a & b & X\setminus\{a,b\} \\
\hline
a & t_a & 1-t_b & O \\
b & 0 & 1 & O \\
X\setminus\{a,b\} & O & O & I
\end{array} \right]\]

We show next that when the rMVA invariant is restricted to $pv\mathcal{B}_X$, it is equivalent to the Gassner representation (and therefore also specializes to the Burau representation when we set all the variables to be equal).  
\begin{proposition}
	On pure v-braids $pv\mathcal{B}_X$, rMVA reduces to the Gassner representation, respectively the Burau representation in the single variable case (letting $t_a=t, \; \forall \; a \in X$).
\end{proposition}
\begin{proof}The image of the space $pv\mathcal{B}_X$ under rMVA has a group structure given by $(\la_1,\mathcal{A}_1)\cdot(\la_2,\mathcal{A}_2)=(\la_1\la_2,\mathcal{A}_1\mathcal{A}_2)$, with matrix multiplication in the second factor. We will first verify that the rMVA invariant, with a small correction, is a group homomorphism. We start by showing that:\begin{equation}\rMVA(B\cdot \sigma_{cd})=-\rMVA(\sigma_{cd})\rMVA(B)\label{eqn:gpmor}\end{equation}
	for any braid $B \in pv\mathcal{B}_{X\cup\{a,b\}}$ and generator $\sigma_{cd}$. Note that the operation of multiplying a pair $(\la,\mathcal{A})$ by a scalar commutes with both metamonoid multiplication, i.e. gluing, and disjoint union as discussed in Section \ref{sub:glue}. So, we will suppress the normalizing factor $1/\sqrt{t_c}$ for $\sigma_{cd}$ in the following computations for the sake of simplicity. The left side of Equation \ref{eqn:gpmor} then becomes:
	\begin{align*}
		&{\small \begin{array}{c|cc}
				t_c & c & d \\
				\hline
				c & -t_c & 0 \\
				d & 1-t_d & -1 
			\end{array}
			\sqcup
			\begin{array}{c|ccc}
				\lambda & a & b & X \\
				\hline
				a & \alpha & \beta & \theta \\
				b & \gamma & \delta & \epsilon \\
				X & \phi & \psi & \Xi 
			\end{array}=\begin{array}{c|ccccc}
			t_c\la & c & d & a & b & X \\
			\hline
			c & -t_c\la  & 0 & 0 & 0 & 0_{|X|} \\
			d & (1-t_d)\la & -\la & 0 & 0 & 0_{|X|} \\
			a & 0 & 0 & t_c\al & t_c\be & t_c\theta \\
			b & 0 & 0 & t_c\ga & t_c\de & t_c\ep \\
			X & 0^{\text{tr}}_{|X|} & 0^{\text{tr}}_{|X|} & t_c\phi & t_c\psi & t_c\Xi
		\end{array}}
	\end{align*}
	
	\begin{align*} & \xrightarrow{m^{ac}_a} \begin{array}{c|cccc}
			t_a\la & a & d & b & X \\
			\hline
			a &  t_a\al & 0 & t_a\beta & t_a\theta \\
			d & (t_d-1)\al & -\la & (t_d-1)\be & (t_d-1)\theta \\
			b & t_a\ga & 0 & t_a\de & t_a\ep \\
			X & t_a\phi & 0^{\text{tr}}_{|X|} & t_a\psi & t_a\Xi
		\end{array} \\
		& \\
		&\xrightarrow{m^{bd}_b} \hspace{-0.1cm}\begin{array}{c|ccc}
		t_a\la & a & b & X \\
		\hline
		a &  t_a\al & t_a\beta & t_a\theta \\
		b & (t_b-1)\al+\ga & (t_b-1)\be+\de & (t_b-1)\theta+\ep \\
		X & t_a\phi & t_a\psi & t_a\Xi
	\end{array}
\end{align*}

Strictly speaking, rather than computing $m^{bd}_d \circ m^{ac}_a((\la_{\si_{cd}},\mathcal{A}_{\si_{cd}})\sqcup(\la_{B},\mathcal{A}_{B}))$, if we think of $\sigma_{cd}$ as a braid on $|X|+2$ strands, then we need to compute $m^{x_{|X|}y_{|X|}}_{x_{|X|}}\circ \hdots \circ m^{x_1y_1}_{x_1}\circ m^{bd}_d \circ m^{ac}_a(e_{y_{|X|}}\circ\hdots\circ e_{y_1}(\la_{\si_{cd}},\mathcal{A}_{\si_{cd}})\sqcup(\la_{B},\mathcal{A}_{B}))$, i.e. include the image of $|X|$ unknotted strands using rMVA. However, the metamonoid axioms tell us that operations on non-overlapping labels commute and $m^{xy}_x \circ e_y=\text{Id}$. So, the computation reduces to the one given above.

On the other hand, to compute the product of the matrices for the two braids on the right side of Equation \ref{eqn:gpmor}, we need them to be of the same size and so we need to write $\sigma_{cd}$ formally as a braid on $|X|+2$ strands. For that purpose we take, on the level of the image of the invariant, the disjoint union of the positive crossing with labels $c$ and $d$, and $|X|$ unknotted strands. As before, we suppress the normalizing factor $1/\sqrt{t_c}$. 

\[{\small \begin{array}{c|cc}
	t_c & c & d \\
	\hline
	c & -t_c & 0 \\
	d & 1-t_d & -1 
	\end{array}\sqcup
	\bigsqcup_{x \in X}\begin{array}{c|c}
	1 & x \\
	\hline
	x & -1
	\end{array}=
	\begin{array}{c|ccc}
	t_c & c & d & X \\
	\hline
	c & -t_c & 0 & 0_{|X|} \\
	d & 1-t_d & -1 & 0_{|X|} \\
	X & 0^{\text{tr}}_{|X|} & 0^{\text{tr}}_{|X|} & -t_cI_{|X| \times |X|}
	\end{array}}\]

Then, the negative of the matrix product (in opposite order) agrees with the matrix from the earlier operation corresponding to stitching together the two braids:

\begin{align*}(\la_{\si_{cd}},\mathcal{A}_{\si_{cd}})\cdot(\la_B,\mathcal{A}_B)=&\left(t_c\la,
	\left[\arraycolsep=3pt\def\arraystretch{1.3}\begin{array}{ccc}
		-t_c & 0 & 0_{|X|} \\
		1-t_d & -1 & 0_{|X|} \\
		0^{\text{tr}}_{|X|} & 0^{\text{tr}}_{|X|} & -t_cI_{|X| \times |X|}
	\end{array}\right]\cdot\left[\arraycolsep=3pt\def\arraystretch{1.3}\begin{array}{ccc}
	\al & \be & \theta \\
	\ga & \de & \ep \\
	\phi & \psi & \Xi 
\end{array}\right]\right) \\
& \\
\xrightarrow{t_c\rightarrow t_a,t_d\rightarrow t_b}&\arraycolsep=3pt\def\arraystretch{1.3}\begin{array}{c|ccc}
	t_a\la & a & b & X \\
	\hline
	a & -t_a\al & -t_a\be & -t_a\theta \\
	b & (1-t_b)\al-\ga & (1-t_b)\be-\de & (1-t_b)\theta-\ep \\
	X & -t_a\phi & -t_a\psi & -t_a\Xi
\end{array}\end{align*}
Both the above product and the matrix from the earlier computation need to be multiplied by the normalizing factor $1/\sqrt{t_a}$ at the end. Thus, taking the negative transpose (or negative inverse) of the image makes rMVA a group homomorphism on pure, virtual braids. Analogous computations can be performed for $\sigma^{-1}_{cd}$. Since the $\sigma_{cd}, \; c,d \in Y:=X\cup\{a,b\}$, generate $pv\mathcal{B}_Y$, we get more generally:
\[\rMVA(B_1\cdot B_2) =-\rMVA(B_2)\cdot\rMVA(B_1) \quad \forall \; B_1,B_2 \in pv\mathcal{B}_Y\]
Furthermore, the image of the generators $\sigma_{cd}$ under rMVA is equivalent to that for the Gassner representation, after taking the negative transpose:
\[\rMVA(\sigma_{cd})=t^{-\frac{1}{2}}_c\begin{array}{c|cc}
t_c & c & d \\
\hline
c & -t_c & 0 \\
d & 1-t_d & -1 
\end{array}\xrightarrow{-(\;\cdot\;)^{\text{tr}}}
t^{-\frac{1}{2}}_c\begin{array}{c|cc}
t_c & c & d \\
\hline
c & t_c & t_d-1 \\
d & 0 & 1 
\end{array}\]
\end{proof}

\begin{remark}
	Since the rMVA invariant satisfies the ``Overcrossings Commute'' relation portrayed in Figures \ref{fig:R1OC} and \ref{fig:ocbc}, it is also a Gassner type invariant for pure welded braids $pw\mathcal{B}_X$.
\end{remark}

	\section{Extensions}
\label{sec:ext}
\subsection{Partial trace} We have mainly discussed the rMVA invariant in the context  of pure tangles, where the gluing operation is well defined, and rMVA is a metamonoid morphism that recovers the tMVA. A brief digression in Section \ref{sec:MVA} discussed the special case of $1$-tangles with possibly some additional closed components. Since the rMVA invariant is obtained from tMVA, which allows gluing of the ends of the same strand and can be evaluated on non-pure tangles albeit in the circuit algebra setting, there is hope that the rMVA invariant can also be extended to non-pure tangles. A step in that direction is the following ``strand closure'' or trace operation induced from the corresponding operation on the target space of tMVA. 
\begin{figure}[H]\arraycolsep=4pt\def\arraystretch{1.2}
	\[\begin{array}{c|cc}
	\la & a & X \\
	\hline
	a & \al & \theta \\
	X & \phi & \Xi
	\end{array}
	\xrightarrow{\quad\tr_a\quad}\arraycolsep=5pt\def\arraystretch{1.8}
	\begin{array}{c|c}
	\la + \al & X \\
	\hline 
	X & \frac{(\la+\al)\Xi - \phi\theta}{\la}	
	\end{array}	\]
	\caption{A partial trace operation on the target space of rMVA.}
	\label{fig:trace}
\end{figure}
\vspace{-0.3cm}
It is only partially defined, as after applying this operation to a pair $(\la,\mathcal{A})$ to get $(\tr_a(\la),\tr_a(\mathcal{A}))$, it is no longer guaranteed that $\tr_a(\la)\neq 0$. So, we cannot apply the metamonoid operations, and we can no longer recover the tMVA using Theorem \ref{thm:deg01}. This operation nonetheless appears promising as it and recovers the MVA for all links with up to $7$ crossings in the Knot Atlas \cite{DBN5}. Hence, if the domain of the trace is understood better, it could allow one to generalize the rMVA invariant to non-pure tangles, i.e. also including closed components.

\subsection{Future directions} Several further questions can be pursued in order to understand the rMVA invariant better and explore possible generalizations. Among them are:
\begin{enumerate}[i)]
	\item The rMVA invariant fits in the setting of metamonoids but is originally obtained through Archibald's tMVA invariant which is a circuit algebra morphism. So, a more general algebraic description of the relation between oriented circuit algebras and metamonoids would aid in understanding the connection better.
	\item To clarify more the connection of rMVA with tMVA and the original multivariable Alexander polynomial on links, it would be beneficial to study further the trace operation in Figure \ref{fig:trace} and the circumstances under which we get $\la=0$ in a pair $(\la,\mathcal{A})$ in the target space of rMVA.
	\item As described in Theorem \ref{thm:iso}, the rMVA is equivalent to a version of Bar-Natan's $Z$ invariant, so the partial trace operation can be considered in that context. The invariant $Z$ itself is a reduction of an invariant of ribbon-knotted copies of $S^1$ and $S^2$ in $\mathbb{R}^4$ mapping to certain free Lie and cyclic words, defined by Bar-Natan in \cite{DBN3}. So, it would be of interest to explore how the trace manifests there.
	\item The Alexander polynomial satisfies several skein relations, most of which have been verified by Archibald in \cite{Arch10} through the tMVA invariant in a more straightforward way than the standard setting, and whose proofs might be further simplified using rMVA.
	\item Finally, since the image of a tangle under the studied invariants is essentially a collection of Laurent polynomials, it would be interesting to seek a categorification of rMVA and $Z$.
\end{enumerate}
	
	\section{Proofs}
\label{sec:proof}

\subsection{Invariance of rMVA}

Similarly to the ``Overcrossings Commute'' move, for the Reidemeister $2$ and $3$ moves we also distinguish between cyclic and braid-like depending on whether the region bounded by the strands can be consistently oriented or not. For each of these relations, we verify that the value of rMVA on the left side matches that on the right side.
\vspace{0.2cm}

\underline{\sl Reidemeister 2 moves.} 

\begin{figure}[H]
	\begin{center}
		\hspace{-0.5cm}\begin{minipage}[c]{0.30\textwidth}
			\def\svgwidth{4cm}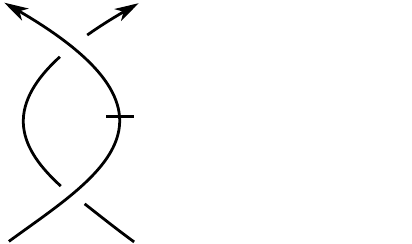
		\end{minipage}
		\hspace{2cm}
		\begin{minipage}[c]{0.30\textwidth}
			\def\svgwidth{4cm}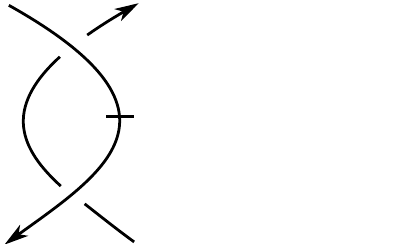
		\end{minipage}
	\end{center}
	\caption{The braid-like and cyclic Reidemeister 2 moves.}
	\label{fig:r2bc}
\end{figure}

For the diagram on the right side of both Reidemeister 2 moves, we have:

\begin{align*}
	&M_R =	\arraycolsep=2pt\def\arraystretch{1.2}\begin{array}{c|cccc}
		& b_1 & b_2 & a_1 & a_2 \\
		\hline
		b_1 & 1 & 0 & -1 & 0 \\
		b_2 & 0 & 1 & 0 & -1
	\end{array} \Rightarrow
	\la_R=1 \; 
	A_R=\begin{array}{c|cc}
		& a_1 & a_2 \\
		\hline
		b_1 & -1 & 0 \\
		b_2 & 0 & -1
	\end{array}  \\
	&\Rightarrow
	\rMVA(D_R)=\prod^2_{s=1}{t^{-\frac{\mu(s)}{2}}_s}\cdot
	\arraycolsep=2pt\def\arraystretch{1.2}\begin{array}{c|c}
		\la_R & \Xin \\
		\hline
		\Xo & A_R
	\end{array}=
	1 \cdot
	\begin{array}{c|cc}
		1 & a_1 & a_2 \\
		\hline
		b_1 & -1 & 0 \\
		b_2 & 0 & -1
	\end{array}=
	\begin{array}{c|cc}
		1 & a_1 & a_2 \\
		\hline
		b_1 & -1 & 0 \\
		b_2 & 0 & -1
	\end{array}
\end{align*}

For the braid-like move, the left diagram $D_{LB}$ gives the following Alexander matrix and rMVA value:
\begin{align*}
	&M_{LB} = \arraycolsep=2pt\def\arraystretch{1.2} \begin{array}{c|ccccc}
		& c & b_1 & b_2 & a_1 & a_2 \\
		\hline
		c & t_1 & 0 & 0 & 1-t_2 & -1 \\
		b_1 & 0 & 1 & 0 & -1 & 0 \\
		b_2 & -t_1 & t_2-1 & 1 & 0 & 0
	\end{array} \Rightarrow
	\la_{LB}=t_1, \; 
	A_{LB}=\begin{array}{c|cc}
		& a_1 & a_2 \\
		\hline
		b_1 & -t_1 & 0 \\
		b_2 & 0 & -t_1
	\end{array} \\
	&\Rightarrow
	\rMVA(D_{LB})=\prod^2_{s=1}{t^{-\frac{\mu(s)}{2}}_s}\cdot
	\arraycolsep=2pt\def\arraystretch{1.2} \begin{array}{c|c}
		\la_{LB} & \Xin \\
		\hline
		\Xo & A_{LB}
	\end{array}=
	t^{-1}_1 \cdot
	\begin{array}{c|cc}
		t_1 & a_1 & a_2 \\
		\hline
		b_1 & -t_1 & 0 \\
		b_2 & 0 & -t_1
	\end{array}=
	\begin{array}{c|cc}
		1 & a_1 & a_2 \\
		\hline
		b_1 & -1 & 0 \\
		b_2 & 0 & -1
	\end{array}
\end{align*}

Similarly, for the cyclic move, the left side $D_{LC}$ gives the following Alexander matrix and rMVA value:
\begin{align*}
	&M_{LC} = \arraycolsep=2pt\def\arraystretch{1.2} \begin{array}{c|ccccc}
		& c & b_1 & b_2 & a_1 & a_2 \\
		\hline
		c & 1 & t_2-1 & 0 & 0 & -t_1 \\
		b_1 & 0 & 1 & 0 & -1 & 0 \\
		b_2 & -1 & 0 & t_1 & 1-t_2 & 0
	\end{array} \Rightarrow
	\la_{LC}=t_1, \; 
	A_{LC}=\begin{array}{c|cc}
		& a_1 & a_2 \\
		\hline
		b_1 & -t_1 & 0 \\
		b_2 & 0 & -t_1
	\end{array} \\
	&\Rightarrow
	\rMVA(D_{LC})=\prod^2_{s=1}{t^{-\frac{\mu(s)}{2}}_s}\cdot
	\arraycolsep=2pt\def\arraystretch{1.2} \begin{array}{c|c}
		\la_{LC} & \Xin \\
		\hline
		\Xo & A_{LB}
	\end{array}=
	t^{-1}_1 \cdot
	\begin{array}{c|cc}
		t_1 & a_1 & a_2 \\
		\hline
		b_1 & -t_1 & 0 \\
		b_2 & 0 & -t_1
	\end{array}=
	\begin{array}{c|cc}
		1 & a_1 & a_2 \\
		\hline
		b_1 & -1 & 0 \\
		b_2 & 0 & -1
	\end{array}
\end{align*}

Next, we consider the braid-like and cyclic Reidemeister 3 moves.
\vspace{0.2cm}

\underline{\sl Reidemeister 3 moves.}

\begin{figure}[H]
	\begin{center}	
		\hspace{-1cm}\begin{minipage}[c]{0.30\textwidth}
			\def\svgwidth{5cm}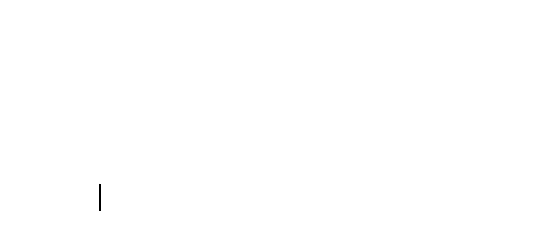
		\end{minipage}
		\hspace{3cm}
		\begin{minipage}[c]{0.30\textwidth}
			\def\svgwidth{5cm}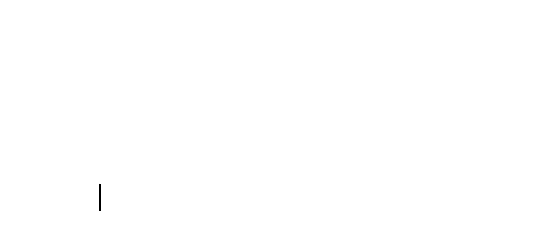
		\end{minipage}
	\end{center}
	\caption{The braid-like and cyclic Reidemeister 3 moves.}
	\label{fig:r3bc}
\end{figure}

Starting with the braid-like move, the diagram on the left side $D_{LB}$ gives the following Alexander matrix and corresponding value of the rMVA:
\begin{align*}
	M_{LB} &=\arraycolsep=4pt\def\arraystretch{1.4} \begin{array}{c|ccccccc}
		& c & b_1 & b_2 & b_3 & a_1 & a_2 & a_3 \\
		\hline
		c & 1 & 0 & 0 & 0 & 0 & -t_3 & t_2-1  \\
		b_1 & 0 & 1 & 0 & t_1-1 & -t_3 & 0 & 0 \\
		b_2 & -1 & 1-t_2 & t_1 & 0 & 0 &0 & 0 \\
		b_3 & 0 & 0 & 0 & 1 & 0 & 0 & -1
	\end{array} \\
	& \\
	\Rightarrow
	\rMVA(D_{LB})&=\prod^3_{s=1}{t^{-\frac{\mu(s)}{2}}_s}\cdot
	\arraycolsep=4pt\def\arraystretch{1.4} \begin{array}{c|c}
		\la_{LB} & \Xin \\
		\hline
		\Xo & A_{LB}
	\end{array} \\
	& \\
	&=
	t^{-\frac{1}{2}}_1t^{-1}_3 \cdot
	\arraycolsep=4pt\def\arraystretch{1.4} \begin{array}{c|ccc}
		t_1 & a_1 & a_2 & a_3\\
		\hline
		b_1 & -t_1t_3 & 0 & t_1(t_1-1) \\
		b_2 & t_3(1-t_2) & -t_3 & t_1(t_2-1) \\
		b_3 & 0 & 0 & -t_1
	\end{array}
\end{align*}

The diagram on the right side of the braid-like Reidemeister 3 move, $D_{RB}$, produces analogously:

\begin{align*}
	M_{RB} &= \arraycolsep=4pt\def\arraystretch{1.4}	\begin{array}{c|ccccccc}
		& c & b_1 & b_2 & b_3 & a_1 & a_2 & a_3 \\
		\hline
		c & t_1 & 0 & 0 & 0 & 1-t_2 & -1 & 0  \\
		b_1 & 0 & 1 & 0 & 0 & -t_3 & 0 & t_1-1 \\
		b_2 & -t_3 & 0 & 1 & t_2-1 & 0 &0 & 0 \\
		b_3 & 0 & 0 & 0 & 1 & 0 & 0 & -1
	\end{array} \\
	& \\
	\Rightarrow
	\rMVA(D_{RB})&=\prod^3_{s=1}{t^{-\frac{\mu(s)}{2}}_s}\cdot
	\arraycolsep=4pt\def\arraystretch{1.4} \begin{array}{c|c}
		\la_{RB} & \Xin \\
		\hline
		\Xo & A_{RB}
	\end{array} \\
	& \\
	&=
	t^{-\frac{1}{2}}_1t^{-1}_3 \cdot
	\arraycolsep=4pt\def\arraystretch{1.4} \begin{array}{c|ccc}
		t_1 & a_1 & a_2 & a_3\\
		\hline
		b_1 & -t_1t_3 & 0 & t_1(t_1-1) \\
		b_2 & t_3(1-t_2) & -t_3 & t_1(t_2-1) \\
		b_3 & 0 & 0 & -t_1
	\end{array}
\end{align*}

Looking next at the second type of Reidemeister 3 move, the cyclic one, the diagram on the left side of the move $D_{LC}$ gives rise to the following Alexander matrix and rMVA value:

\begin{align*}
	M_{LC} &= \arraycolsep=4pt\def\arraystretch{1.4} \begin{array}{c|ccccccc}
		& c & b_1 & b_2 & b_3 & a_1 & a_2 & a_3 \\
		\hline
		c & 1 & t_2-1 & 0 & 0 & 0 & -t_1 & 0  \\
		b_1 & 0 & 1 & 0 & t_1-1 & -t_3 & 0 & 0 \\
		b_2 & -1 & 0 & t_3 & 0 & 0 & 0 & 1-t_2 \\
		b_3 & 0 & 0 & 0 & 1 & 0 & 0 & -1
	\end{array} \\
	& \\
	\Rightarrow
	\rMVA(D_{LC})&=\prod^3_{s=1}{t^{-\frac{\mu(s)}{2}}_s}\cdot
	\arraycolsep=4pt\def\arraystretch{1.4} \begin{array}{c|c}
		\la_{LC} & \Xin \\
		\hline
		\Xo & A_{LC}
	\end{array} \\
	& \\
	&=
	t^{-\frac{1}{2}}_1t^{-1}_3 \cdot
	\arraycolsep=4pt\def\arraystretch{1.4} \begin{array}{c|ccc}
		t_3 & a_1 & a_2 & a_3\\
		\hline
		b_1 & -t^2_3 & 0 & t_3(t_1-1) \\
		b_2 & t_3(t_2-1) & -t_1 & t_1(1-t_2) \\
		b_3 & 0 & 0 & -t_3
	\end{array}
\end{align*}

For the diagram on the right side of the same move, $D_{RC}$, we get:

\begin{align*}
	M_{RC} &= \arraycolsep=4pt\def\arraystretch{1.4} \begin{array}{c|ccccccc}
		& c & b_1 & b_2 & b_3 & a_1 & a_2 & a_3 \\
		\hline
		c & t_3 & 0 & 0 & 1-t_2 & 0 & -1 & 0  \\
		b_1 & 0 & 1 & 0 & 0 & -t_3 & 0 & t_1-1 \\
		b_2 & -t_1 & 0 & 1 & 0 & t_2-1 & 0 & 0 \\
		b_3 & 0 & 0 & 0 & 1 & 0 & 0 & -1
	\end{array} \\
	& \\
	\Rightarrow
	\rMVA(D_{RC})&=\prod^3_{s=1}{t^{-\frac{\mu(s)}{2}}_s}\cdot
	\arraycolsep=4pt\def\arraystretch{1.4} \begin{array}{c|c}
		\la_{RC} & \Xin \\
		\hline
		\Xo & A_{RC}
	\end{array} \\
	& \\
	&=
	t^{-\frac{1}{2}}_1t^{-1}_3 \cdot
	\arraycolsep=4pt\def\arraystretch{1.4} \begin{array}{c|ccc}
		t_3 & a_1 & a_2 & a_3\\
		\hline
		b_1 & -t^2_3 & 0 & t_3(t_1-1) \\
		b_2 & t_3(t_2-1) & -t_1 & t_1(1-t_2) \\
		b_3 & 0 & 0 & -t_3
	\end{array}
\end{align*}

\subsection{Proof of Theorem \ref{thm:deg01}}

Theorem \ref{thm:deg01} follows from the proof of the generalized Cramer's rule by Gong, Aldeen and Elsner \cite{GAE}, which we reproduce here in our setting.
Throughout, we will use the notation:

\[
\begin{array}{lll}
A^{j_1,\hdots,j_k}_{i_1,\hdots,i_k}  & \hspace{-0.2cm}= &\text{the $k \times k$ submatrix of $A$ with} \\
& & \text{ columns $j_1,\hdots,j_k$ and rows $i_1,\hdots,i_k$} \\
&& \\
A^{p_1,\hdots,p_k}_{B(q_1,\hdots,q_k)}  & \hspace{-0.2cm}= &\text{the matrix $A$ with column $p_s$} \\
& & \text{ replaced by column $q_s$ of $B$, $s=1,\hdots,k$}
\end{array}
\]

\vspace{0.2cm}

Suppose $T \in pv\mathcal{T}_X$ is a pure regular v-tangle with $m$ internal arcs and $|X|=n$, with $(m+n) \times (m+2n)$ Alexander matrix $M(D_T)$ for a diagram $D_T$ of $T$. For convenience, we will distinguish $\Xin \cong \Xo \cong X$, the labels of the incoming and outgoing arcs of $T$. As discussed earlier, the degree $0$ and $1$ components for $T$ using the tMVA invariant and the Hodge $\ast$ operator are:
\[\la=\det{M(D_T)^{1,\hdots,n;\emptyset}} \quad \quad  \mathcal{A}_{i,j}=(-1)^{n-i}\det{M(D_T)^{1,\hdots,\hat{i},\hdots,n;j}},\; 1\leq i,j \leq n \vspace{0.2cm}\]
Take the $(m+n) \times (m+n)$ matrix $D=M(D_T)^{1,\hdots,n;\emptyset}$ ($\det{D}=\lambda$) and the $(m+n) \times n$ matrix $N$, $N_{i,j}=M(D_T)_{i,(m+n)+j}$, namely the matrix made up of the columns of the Alexander matrix labelled by $\Xin$. 

Now, consider the $(m+n) \times n$ matrix $H$ defined by $H_{i,j}=\frac{\det{D^i_{N(j)}}}{\det{D}}$. Note that $H_{m+i,j}=\frac{\mathcal{A}_{i,j}}{\det{D}}$ for $1 \leq i \leq n$.

\begin{lemma}\label{lm:mprod} Let $I$ denote the $(m+n) \times (m+n)$ identity matrix. Then:
	\begin{align}
		D \cdot I^{i_1+m,\hdots,i_k+m}_{H(j_1,\hdots,j_k)} &= D^{i_1+m,\hdots,i_k+m}_{N(j_1,\hdots,j_k)}
		\label{eq:mprod}
	\end{align}
\end{lemma}

\begin{proof}
	For any matrix $C$, let $C_{\ast,j}$ denote its $j$th column vector. Then the equation in the lemma can be visualized as:	
	\[\begin{array}{c}
	\kbordermatrix{ & &\text{int} 	 & 	 & 				&  & \Xo  \cr
		& & \phantom{1}  & \phantom{\vline height 3ex} & \vrule &  &  \cr
		\rotatebox{90}{\tiny int} & & \phantom{0}	 & 	& \vrule &  &  \cr	 
		& & \phantom{0}	 & 	& \vrule &  &  \cr \cline{2-7}
		& & \phantom{0} & 	 & \vrule &  &  \cr 
		\rotatebox{90}{$\scriptstyle \Xo$} & & \phantom{0} & \phantom{\vline height 3ex} & \vrule &  & }
	\times\hspace{-0.25cm}
	\kbordermatrix{ & & i_1+m & & i_k+m & \cr
		& 1 & \vline height 3ex	& & \vline height 3ex & 0 \cr
		& 0 &	&   &  & 0	\cr	
		& \vdots & H_{\ast,j_1}	& \cdots & H_{\ast,j_k}	 & \vdots \cr 
		& 0 &	&  &  & 0	\cr								
		& 0 & \vline height 3ex	&  & \vline height 3ex & 1 \cr }
	= \\
	=\kbordermatrix{& &\text{int}   & 	& & i_1+m	& & i_k+m\cr
		& & \phantom{1} & \phantom{\vline height 3ex} & \vrule & \vline height 3ex &  & \vline height 3ex \cr
		& & \phantom{0}	& &\vrule &  &  \cr	
		& & \phantom{\vdots} &  & \vrule & N_{\ast,j_1} & \cdots & N_{\ast,j_k}\cr \cline{2-8}
		& & \phantom{0} & & \vrule & & 	& \cr 
		& & \phantom{0} & \phantom{\vline height 3ex} & \vrule & \vline height 3ex 	& & \vline height 3ex}\end{array}\]
	Let $B=I^{i_1+m,\hdots,i_k+m}_{H(j_1,\hdots,j_k)}$. Then the equation follows directly from the classical Cramer's rule, from which we can conclude that:
	\[D \cdot B_{\ast,s} = \begin{cases} 
	D \cdot H_{\ast,j_p}=D \cdot \frac{\det{D^{\ast}_{N(j_p)}}}{\det{D}} = N_{\ast,j_p} & s=m+i_p, \; (p=1,\hdots,k) \\
	D \cdot I_{\ast,s}=D_{\ast,s} & \text{otherwise.} \\
	\end{cases}\]
\end{proof}

We take determinants of Equation \ref{eq:mprod} in Lemma \ref{lm:mprod}, and use:

\begin{lemma} For any $1 \leq i_1,\hdots,i_k,j_1,\hdots,j_k \leq n$:
	\[
	\det{I^{i_1+m,\hdots,i_k+m}_{H(j_1,\hdots,j_k)}}=\frac{1}{\lambda^k}\det{\mathcal{A}^{j_1,\hdots,j_k}_{i_1,\hdots,i_k}}
	\]
\end{lemma}

\begin{proof}
	Let $B=I^{i_1+m,\hdots,i_k+m}_{H(j_1,\hdots,j_k)}$. Then, using $\delta_{x,y}$ to denote the Kronecker delta ($\delta_{x,y}=1$ if $x=y$ and is zero otherwise):
	
	\begin{align*}
		\det{B} &= \hspace{-0.2cm} \sum_{\sigma \in S_{n+m}} (-1)^{\sigma} b_{\sigma(1),1} \hdots b_{\sigma(i_1+m),i_1+m} \hdots b_{\sigma(i_k+m),i_k+m} \hdots b_{\sigma(n+m),n+m} \\
		&= \hspace{-0.2cm} \sum_{\sigma \in S_{n+m}} (-1)^{\sigma} \delta_{\sigma(1),1} \hdots b_{\sigma(i_1+m),i_1+m} \hdots b_{\sigma(i_k+m),i_k+m} \hdots \delta_{\sigma(n+m),n+m} \\
		&= \sum_{\sigma \in S_k} (-1)^{\sigma} b_{\sigma(i_1+m),i_1+m} b_{\sigma(i_2+m),i_2+m} \hdots b_{\sigma(i_k+m),i_k+m} \\
		&= \det{H^{j_1,\hdots,j_k}_{i_1+m,\hdots,i_k+m}} =\frac{1}{\lambda^k} \det{\mathcal{A}^{j_1,\hdots,j_k}_{i_1,\hdots,i_k}}
	\end{align*}
\end{proof}

The determinant version of Equation \ref{eq:mprod} then becomes:
\begin{align*}
	\la \frac{1}{\la^k}\det{\mathcal{A}^{j_1,\hdots,j_k}_{i_1,\hdots,i_k}} &= \det{D^{i_1+m,\hdots,i_k+m}_{N(j_1,\hdots,j_k)}} \\
&= (-1)^{nk-\frac{(k-1)k}{2}-\sum^k_{p=1}{i_p}}\det{M(D_T)^{\{1,\hdots,n\} \setminus \{i_1,\hdots,i_k\};j_1,\hdots,j_k}}
\end{align*}

The coefficients of the tMVA invariant are thus obtained from the pair $(\lambda,\mathcal{A})$ whenever $\la \neq 0$, using the formula:
\[
\boxed{\det{M(D_T)^{\{1,\hdots,n\} \setminus \{i_1,\hdots,i_k\};j_1,\hdots,j_k}}=(-1)^{nk-\frac{(k-1)k}{2}-\sum^k_{p=1}{i_p}}\frac{\det{\mathcal{A}^{j_1,\hdots,j_k}_{i_1,\hdots,i_k}}}{\lambda^{k-1}}}
\]

\subsection{Proof of Theorem \ref{thm:deg01maps}} Suppose $T$ is a pure regular virtual tangle. Then, using the earlier notation: 
\begin{align*}&\tMVA(T)=\prod_k{t^{-\frac{\mu(k)}{2}}_k} w \otimes \sum^n_{k=0}\sum_{\overline{i}, \underline{j}} {\det{M(D_T)^{\overline{i};\underline{j}}}\xin_{\underline{j}} \otimes  \xo_{\overline{i}}} \\
	&\Rightarrow p_2(\tMVA(T))=\sum^n_{k=0}\sum_{\overline{i},\underline{j}} {\det{M(D_T)^{\overline{i};\underline{j}}}\xin_{\underline{j}} \otimes  \xo_{\overline{i}}}
\end{align*}
To shorten the notation, we denote as before $\left[n\right]=\{1,\hdots,n\}$, $\overline{i}=\{i_1<\hdots<i_{n-k}\},\;\underline{j}=\{j_1 <\hdots< j_k\}$, $x_{\overline{i}}=x_{i_1} \wedge \hdots \wedge x_{i_{n-k}}$, and $y_{\underline{j}}=y_{j_1}\wedge\hdots\wedge y_{j_k}$, as well as $x_{\underline{s}}(y_{\underline{t}})=x_{s_1}(y_{t_1})\cdot\hdots\cdot x_{s_k}(y_{t_k})$, as $\Xin$ and $\Xo$ are self-dual. Following the diagram in Equation \ref{eqn:Hodge}, for $p_2(\tMVA(T))$, we get:
\begin{align*}&p_2(\tMVA(T))\xrightarrow{\ast_w} \sum^n_{k=0}\sum_{\overline{i},\underline{j}}{(-1)^{\ast_w}{\det M(D_T)^{\overline{i};\underline{j}}} \xin_{\underline{j}} \otimes \xo_{\left[n\right]\setminus\overline{i}} }\rightarrow \\
	&\rightarrow \left(\phi_k:\xin_{\underline{p}}\mapsto \sum_{\overline{i},\underline{j}}{(-1)^{\ast_w}{\det M(D_T)^{\overline{i};\underline{j}}}}\left[\sum_{\si \in S_k}\text{sgn}(\si)\xin_{\si(\underline{j})}(\xin_{\underline{p}}) \right]\xo_{\left[n\right]\setminus\overline{i}}\right)^n_{k=0}\hspace{-0.5cm}= \\
	&=\left(\phi_k:\xin_{\underline{p}}\mapsto \sum_{\overline{i}}{(-1)^{\ast_w}{\det M(D_T)^{\overline{i};\underline{p}}}}\xo_{\left[n\right]\setminus\overline{i}}\right)^n_{k=0} 
\end{align*}
Here, the image under $\ast_w$ belongs to $\bigoplus^n_{k=0} \Lambda^k(\Xin) \otimes \Lambda^k(\Xo)$, and the final image belongs to $\bigoplus^n_{k=0} \text{Hom}(\Lambda^k(\Xin), \Lambda^k(\Xo))$. In particular, the degree $0$ and $1$ components are:
\begin{align*}
	\phi_0(1)&=\det{M(D_T)^{\left[n\right];\emptyset}}=\la \\
	\phi(\xin_j):=\phi_1(\xin_j)&=\sum^n_{i=1}(-1)^{n-i}\det{ M(D_T)^{\left[n\right]\setminus\{i\};j}}\xo_i
\end{align*}
Then, taking $\psi=\phi/\la$, we have for any $0 \leq k \leq n$, $1 \leq j_1 < \hdots <j_k \leq n$:
\begin{align*}&\Lambda(\psi)(\xin_{j_1} \wedge \hdots \wedge \xin_{j_k}) =\psi(\xin_{j_1}) \wedge \hdots \wedge \psi(\xin_{j_k}) \\
	& \left[\psi(\xin_{j_s})=\sum^n_{i_s=1}(-1)^{n-i_s}\frac{\det{ M(D_T)^{\left[n\right]\setminus\{i_s\};j_s}}}{\la}\xo_{i_s} \right] \\
	&= \left(\sum^n_{i_1=1}\frac{\mathcal{A}_{i_1,j_1}}{\la}\xo_{i_1} \right)\wedge \hdots \wedge \left(\sum^n_{i_k=1}\frac{\mathcal{A}_{i_k,j_k}}{\la}\xo_{i_k}\right) \\
	& \\
	&=\frac{1}{\la^k}\sum_{\underline{i}}\sum_{\si \in S_k}\text{sgn}(\si)\prod^k_{s=1}\mathcal{A}_{\si(i_s),j_s}\xo_{\underline{i}} =\sum_{\underline{i}}\frac{\det{\mathcal{A}^{\underline{j}}_{\underline{i}}}}{\la^k}\xo_{\underline{i}} \\
	& \\
	&\stackrel{Thm \ref{thm:deg01}}{=}\frac{1}{\la}\sum_{\underline{i}}(-1)^{nk-\sum^k_{p=1}{i_p}-(k-1)k/2}\det{M(D_T)^{\{1,\hdots,n\} \setminus \underline{i};\underline{j}}}\xo_{\underline{i}} \\
	& \\
	&=\frac{1}{\la}\phi_k(\xin_{j_1} \wedge \hdots \wedge \xin_{j_k})
\end{align*}
Therefore, $\la\Lambda(\phi/\la)=(\phi_k)^n_{k=0}$, concluding the proof.

	\end{document}